\documentclass[12pt,oneside,reqno]{amsart}
\usepackage{graphicx}
\usepackage{mathrsfs}
\usepackage{stmaryrd}
\usepackage{amsfonts}
\usepackage{cite}
\usepackage{enumerate,amsmath,amssymb,amsthm}
\usepackage{booktabs} 
\usepackage{diagbox} 
\usepackage{xcolor}

\newcommand{\dif}{\mathrm{d}}

\newcommand{\be}{\begin{eqnarray}}
\newcommand{\ee}{\end{eqnarray}}
\newcommand{\ce}{\begin{eqnarray*}}
\newcommand{\de}{\end{eqnarray*}}
\newtheorem{theorem}{Theorem}[section]
\newtheorem{lemma}[theorem]{Lemma}
\newtheorem{remark}[theorem]{Remark}
\newtheorem{definition}[theorem]{Definition}
\newtheorem{proposition}[theorem]{Proposition}
\newtheorem{Examples}[theorem]{Examples}
\newtheorem{corollary}[theorem]{Corollary}
\newtheorem{condition}[theorem]{Condition}
\def\e{\varepsilon}
\def\t{\theta}

\def\d{\delta}

\def\g{\gamma}
\def\s{\sigma}

\def\[{{\Big[}}
\def\]{{\Big]}}
\def\<{{\langle}}
\def\>{{\rangle}}
\def\({{\Big(}}
\def\){{\Big)}}

\def\no{\nonumber}
\def\bt{\begin{theorem}}
	\def\et{\end{theorem}}
\def\bl{\begin{lemma}}
	\def\el{\end{lemma}}
\def\br{\begin{remark}}
	\def\er{\end{remark}}
\def\bx{\begin{Examples}}
	\def\ex{\end{Examples}}
\def\bd{\begin{definition}}
	\def\ed{\end{definition}}
\def\bp{\begin{proposition}}
	\def\ep{\end{proposition}}
\def\bc{\begin{corollary}}
	\def\ec{\end{corollary}}
\def\bco{\begin{condition}}
	\def\eco{\end{condition}}

\def\cC{{\mathcal C}}
\def\cD{{\mathcal D}}

\def\cG{{\mathcal G}}

\def\cS{{\mathcal S}}

\def\mE{{\mathbb E}}

\def\mN{{\mathbb N}}

\def\mP{{\mathbb P}}

\def\mR{{\mathbb R}}
\def\mS{{\mathbb S}}

\def\mW{{\mathbb W}}

\def\sA{{\mathscr A}}
\def\sB{{\mathscr B}}
\def\sC{{\mathscr C}}

\def\sF{{\mathscr F}}

\def\sL{{\mathscr L}}

\def\sP{{\mathscr P}}

\def\sV{{\mathscr V}}

\def\geq{\geqslant}
\def\leq{\leqslant}

\begin{document}

\allowdisplaybreaks

\title{Small noise asymptotic behaviors for path-dependent multivalued McKean-Vlasov stochastic differential equations}

\author{Ying Ma and Huijie Qiao}



\thanks{This work was supported by NSF of China (No.12071071) and the Jiangsu Provincial Scientific Research Center of Applied Mathematics (No. BK20233002).}

\thanks{Corresponding author: Huijie Qiao, hjqiaogean@seu.edu.cn}


\date{}

\dedicatory{Department of Mathematics,
Southeast University,\\
Nanjing, Jiangsu 211189, P.R.China}

\begin{abstract}
This paper investigates the asymptotic behavior of path-dependent multivalued McKean-Vlasov stochastic differential equations perturbed by small noise. Specifically, we first establish a large deviation principle for such equations under non-Lipschitz coefficients by the weak convergence approach. Subsequently, we introduce an auxiliary equation and apply it to derive the moderate deviation principle. Finally, we construct another auxiliary equation and a limit equation, and prove the central limit theorem.
\end{abstract}

\maketitle \rm

\section{Introduction}

Fix $T>0$ and consider the following multivalued stochastic differential equation (SDE) on $\mR^d$:
\be\left\{\begin{array}{ll}
	\dif X(t)\in -A(X(t))\dif t+b(X_t,\sL_{X_t})\dif t+\s(X_t,\sL_{X_t})\dif W(t), 0\leq t\leq T,\\
	X_0=\xi \in \sC,
\end{array}
\right.
\label{0eq}
\ee
where $A: \mathbb{R}^d\to 2^{\mathbb{R}^d}$ is a maximal monotone operator (See Subsection \ref{maxi}), $\{W(t),t\geq0\}$ is a $m$-dimensional standard Brownian motion defined on a complete filtered probability space $(\Omega,\sF,\{\sF_t\}_{t\geq 0},\mP)$, $\sC$ stands for a continuous function space and $\sP_{2}^{\sC}$ denotes the class of probability measures on $\sC$ which have finite second-order moments (See Subsection \ref{nota}). $b: \sC\times\sP_{2}^{\sC}\to \mR^d$, $ \sigma: \sC\times\sP_{2}^{\sC}\to \mR^{d\times m}$ are Borel measurable and $\xi$ is deterministic. Set $X_t(\theta):=X(t+\theta), \theta\in[-r_0,0]$, and $X_{\cdot}$ is called the segment functional of $X(\cdot)$. And $\sL_{X_t}$ denotes the distribution of $X_t$. Since the coefficients of Eq.(\ref{0eq}) depend on not only past paths but also distributions, we call it a path-dependent multivalued McKean-Vlasov SDE.

Equations like (\ref{0eq}) constitute a natural generalization of McKean-Vlasov SDEs, also referred to in the literature as distribution-dependent SDEs or mean-field SDEs. When the coefficients are independent of the solution's past trajectory, Eq.(\ref{0eq}) reduces to a multivalued McKean-Vlasov SDE; corresponding results about well-posedness (\cite{chi, qg}), ergodicity (\cite{q2}), and large deviation principles (\cite{clqz, flqz}) have been established in prior work. Further, in the case where the multivalued operator $A$ vanishes identically (i.e., $A \equiv 0$), Eq.(\ref{0eq}) specializes to a path-dependent McKean-Vlasov SDE. Recent advances in this setting include rigorous treatments of well-posedness, regularity, and asymptotic behavior (\cite{brw, gs, huangx1, hrockw, hy, szw, zhaox}).

In \cite{mq} we have proved the well-posedness for Eq.(\ref{0eq}) under non-Lipschitz conditions. So, in this paper we continue to study the asymptotic behaviors of small perturbation for path-dependent multivalued McKean-Vlasov SDEs. Concretely speaking, we consider the following path-dependent multivalued McKean-Vlasov SDE:
\be\left\{\begin{array}{ll}
	\dif X^\e(t)\in -A(X^\e(t))\dif t+b(X^\e_t,\sL_{X^\e_t})\dif t+\sqrt{\e}\s(X^\e_t,\sL_{X^\e_t})\dif W(t), \quad t\in(0,T],\\
	X^\e(t)=\xi(t) \in \overline{\cD(A)}, \quad t\in[-r_0,0].
\end{array}
\right.
\label{eq0}
\ee
Assume that $(X^\e_\cdot,K^\e(\cdot))$ is the solution of Eq.(\ref{eq0}) and $(X^0_{\cdot},K^0({\cdot}))$ is the solution of the following deterministic multivalued differential equation:
\ce
\left\{\begin{array}{ll}
	\dif X^0(t)\in -A(X^0(t))\dif t+b(X^0_t,\d_{X^0_t})\dif t, \quad t\in(0,T],\\
	X^0(t)=\xi(t) \in \overline{\cD(A)}, \quad t\in[-r_0,0],
\end{array}
\right.
\de
where $\d_{X^0_t}$ is the Dirac measure at $X^0_t$. Then we intend to study the asymptotic behavior of the trajectory
$$
\frac{X^\e(t)-X^0(t)}{a(\e)}, \quad t\in[-r_0,T],
$$
as $\e\to0$ with different convergence rates $a(\e)$. Two cases of $a(\e)=1$ and $a(\e)=\sqrt \e$ correspond to the large deviation principle (LDP) and the central limit theorem (CLT), respectively. And when $a(\e)$ satisfies 
$$
a(\e)\to0, \quad \frac{\e}{a^2(\e)}\to0\quad \text{as}~ \e\to0,
$$
it corresponds to the moderate deviation principle (MDP). 

It is well known that LDPs and MDPs quantify the asymptotic probabilities of rare events via an explicit rate function. Within the framework of SDEs, the large deviation theory characterizes a deterministic limiting trajectory-the so-called ``typical path"-around which the solution process concentrates with exponentially high probability; moreover, it furnishes precise exponential decay estimates for the probabilities of deviations from this path. Among the methodological approaches to establishing LDPs, the weak convergence method has emerged as a powerful and widely applicable technique (See \cite{bd1,bd2,bd3} and the references therein).

LDPs, MDPs  and CLTs have been extensively investigated in the literature. For multivalued SDEs, Ren, Xu, and Zhang \cite{ren1} established an LDP, while Zhang \cite{zhangh} derived the corresponding MDP-both via the weak convergence method. In the context of McKean-Vlasov SDEs, Dos Reis et al. \cite{dos} and Herrmann et al. \cite{herrmann} independently obtained LDPs under distinct assumptions; a key technical step in their analyses involves replacing the distribution $\sL_{X^\varepsilon(t)}$ appearing in the coefficients with the Dirac measure $\delta_{X^0(t)}$, thereby decoupling the distribution dependence and reducing the problem to a classical setting. Suo and Yuan \cite{sy} extended this decoupling strategy to derive the MDP for McKean-Vlasov SDEs and further established the associated CLT. Building on these advances, Liu et al. \cite{liu} developed a unified analytical framework that systematically facilitates the application of the weak convergence method to derive both LDPs and MDPs for a broad class of mean-field stochastic systems. Leveraging this framework, Gu and Song \cite{gs} obtained LDPs and MDPs for path-dependent McKean-Vlasov SDEs. Subsequently, Shen, Zhou, and Wu \cite{szw} generalized the methodology to establish analogous large- and moderate-deviation results for path-dependent McKean-Vlasov SDEs driven by mixed fractional Brownian motion. Using the weak convergence method, Fan, Yu, and Yuan \cite{fyy} investigated LDPs and MDPs for McKean-Vlasov SDEs driven by fractional Brownian motion and also proved the corresponding CLT. Finally, for multivalued McKean-Vlasov SDEs, the second author and collaborators \cite{flqz} directly applied the weak convergence method to obtain the LDP, adapted the argument from \cite{sy} to derive the MDP, and established the CLT.

In this paper, we first establish the LDP for Eq.(\ref{0eq}) under non-Lipschitz coefficients via the weak convergence method (\cite{liu}). Our main result (Theorem \ref{ldpth}) generalizes Theorem 3.8 in \cite{flqz} and Theorem 2.6 in \cite{gs}, relaxing their Lipschitz-type assumptions to a more general framework. Subsequently, to address the MDP, we confront a fundamental technical challenge: because the maximal monotone operator $A$ is inherently nonlinear, the standard linearization procedure-i.e., subtracting the deterministic limit $X^0(t)$ from the perturbed solution $X^\varepsilon(t)$-fails to yield a well-defined equation for the scaled deviation $\frac{X^\varepsilon(t)-X^0(t)}{a(\varepsilon)}$. To overcome this obstacle, we construct one auxiliary stochastic equation (Eq.(\ref{eqmdp})) whose solution has the same form as the normalized deviation process; leveraging this construction, we directly derive the MDP without recourse to the decoupling strategy (\cite{flqz}). Notably, our proof relies on fourth-order moment estimates of the deviation process, circumventing the use of stopping times employed in prior works (\cite{gs, liu, szw}). This choice avoids the distribution problem associated with stopping-time truncations (\cite{wangf1}). Finally, again because of nonlinearity for $A$, we construct another auxiliary stochastic equation (Eq.(\ref{eqcltdistb})) and the corresponding limit equation (Eq.(\ref{eqcltlimit})), and prove the CLT. 

The remainder of this paper is organized as follows. Section \ref{sec:notations} collects essential notation and auxiliary results, including key lemmas required in subsequent analysis. In Section \ref{sec:ldp}, we first state a general criterion for LDPs-adapted from the weak convergence framework-and then apply it to establish the LDP for path-dependent multivalued McKean-Vlasov SDEs. Section \ref{sec:mdp} and \ref{sec:clt} are devoted to the MDP and CLT, respectively.

The following convention will be used throughout the paper: $C$ with or without indices will denote different positive constants whose values may change from one place to another.

\section{Preliminary}\label{sec:notations}
In this section, we recall some notations, concepts and some useful lemmas.

\subsection{Notations}\label{nota}
We shall use $|\cdot|$ and $\|\cdot\|$ for norms of vectors and matrices, respectively. Let $\<\cdot,\cdot\>$ denote the scalar product in $\mathbb{R}^d$. $U^{*}$ denotes the transpose of the matrix $U$.

Fix $r_0>0$. For any closed subset ${\rm D}\subset\mR^d$, let $\sC:=C([-r_0,0], {\rm D})$ be the collection of all the continuous functions from $[-r_0,0]$ to ${\rm D}$ with the uniform norm $\|\xi\|_\infty:=\sup\limits_{\theta\in[-r_0,0]}|\xi(\theta)|$. For $\gamma(\cdot)\in C([-r_0,\infty), \rm D)$, the segment functional $\gamma_\cdot\in C(\mathbb{R}_{+}, \sC)$ is defined as $\gamma_t(\theta):=\gamma(t+\theta),\ \theta\in[-r_0,0],t\geq 0.$

Let $\sP_{2}^{\sC}$ be the class of probability measures on $\sC$ which have finite second-order moments, that is 
$$
\mu(\|\cdot\|_{\infty}^{2}):=\int_{\sC}\|\xi\|_{\infty}^2\mu(\dif \xi)<\infty, \quad \mbox{for}~ \mu\in\sP_{2}^{\sC}.
$$
Then $\sP_{2}^{\sC}$ is a polish space under the Wasserstein distance
$$
\mW_2\left(\mu,\nu\right):=\inf_{\pi\in\cC(\mu,\nu)}\left(\int_{\sC\times\sC}\|\xi-\eta\|_\infty^2\pi(\dif\xi,\dif\eta)\right)^{\frac{1}{2}},
$$
where $\cC(\mu,\nu)$ is the class of couplings for $\mu$ and $\nu$. Moreover, if $\mu$ and $\nu$ are distributions of $X,Y\in L^2(\Omega, \sF, \mP; \sC)$ respectively, it holds that 
\ce
\mathbb{W}_2(\mu,\nu)^2\leq \mE\|X-Y\|^2_{\infty},
\de
where $\mE$ denotes the expectation with respect to $\mP$.

\subsection{Maximal monotone operators}\label{maxi}
For a multivalued operator $A:\mathbb{R}^d\to 2^{\mathbb{R}^d}$, let
\ce
\mathcal{D}(A)&:=&\{x\in\mathbb{R}^d:A(x)\neq\emptyset\},\\
\mathrm{Gr}(A)&:=&\{(x,y)\in\mathbb{R}^{2d}:x\in\mathcal{D}(A),y\in A(x)\}.
\de
We say that $A$ is monotone if for any $(x_{1},y_{1}),(x_{2},y_{2})\in \mathrm{Gr}(A)$, 
$$
\< x_{1}-x_{2},y_{1}-y_{2}\>\geq 0;
$$
$A$ is maximal monotone if  
$$
(x_1,y_1)\in \mathrm{Gr}(A) \Longleftrightarrow\<x_1-x_2,y_1-y_2\>\geq 0,\forall(x_2,y_2)\in \mathrm{Gr}(A).
$$

For any $T>0$, let $\sV_0$ be the set of all continuous functions $K:[0,T]\to \mathbb{R}^d$ with finite variations and $K(0)=0$. And we denote the variation of $K$ on $[0,s]$ for some $s\in[0,T]$ as $|K|_0^{s}$. Set
\ce	
&&\sA:=\{\left(X,K\right):X\in C([0,T], \overline{\cD(A)}),K\in\sV_{0},\\ 
&&\qquad\quad \text{and}~\< X(t)-x,\dif K(t)-y\dif t\>\geq 0,\forall(x,y)\in\mathrm{Gr}(A)\}.
\de

As for $\sA$, the following properties hold (cf. \cite{Cepa98,zhangx2}).

\bl\label{L1}
For $X\in C([0,T], \overline{\cD(A)})$ and $K\in\sV_0$, the following statements are equivalent:
\begin{enumerate}[(i)]
	\item $(X,K)\in\sA$;
	\item For any $x,y\in C([0,T], \mathbb{R}^{d})$ with $(x(t),y(t))\in\mathrm{Gr}(A)$, it holds that 
	$$
	\< X(t)-x(t),\dif K(t)-y(t)\dif t\>\geq 0;
	$$
	\item For any $(X',K')\in\sA$, it holds that 
	$$
	\< X(t)-X'(t),\dif K(t)-\dif K'(t)\>\geq 0.
	$$
\end{enumerate}
\el

\bl\label{L2}
Assume that $\mathrm{Int}(\cD(A))\not=\emptyset$, then for any $a\in\mathrm{Int}(\cD(A))$, there exists constants $\gamma_1>0$ and $\gamma_2,\gamma_3\geq0$ such that for any $(X,K)\in\sA$ and $0\leq s<t\leq T$,
$$
\int_s^t\< X(r)-a,\dif K(r)\>\geq\gamma_1|K|_s^t-\gamma_2\int_s^t\left|X(r)-a\right|\dif r-\gamma_3(t-s).
$$	
\el

\subsection{Two useful lemmas}

Two following lemmas are from \cite[Lemma 144]{sr} and \cite{bihari}, respectively.

\bl\label{2k}
If $\kappa_i: [0, \infty)\to [0, \infty)$ is strictly increasing, concave and continuous and satisfies $\int_{0^+}\frac{1}{\kappa_i(v)}\dif v=\infty$ and $\kappa_i(0)=0$, $i=1,2$, then $\kappa_1(v)+\kappa_2(v)$ is still strictly increasing, concave and continuous and satisfies 
$$
\int_{0^+}\frac{1}{\kappa_1(v)+\kappa_2(v)}\dif v=\infty, \quad \kappa_1(0)+\kappa_2(0)=0.
$$
\el

\bl[The Bihari inequality]\label{bihari}
Let $v(t), F(t)$ be positive continuous functions on $[0,T]$ and $c\geq 0$. Let $\kappa: \mR_+\to \mR_+$ be a increasing continuous function. If
$$
v(t)\leq c+\int_{0}^{t}F(s)\kappa(v(s))\dif s, \quad \forall 0\leq t\leq T,
$$
then
$$
v(t)\leq \Psi^{-1}\left(\Psi(c)+\int_{0}^{t}F(s)\dif s\right)
$$
holds for all $t\in[0,T]$ with
$$
\Psi(c)+\int_{0}^{t}F(s)\dif s\in D(\Psi^{-1}),
$$
where $\Psi(t)=\int_{\vartheta}^{t}\frac{1}{\kappa(v)}\dif v$ for $\vartheta>0, t\geq 0$, $\Psi^{-1}$ is the inverse function of $\Psi$ on $[0,T]$ and $D(\Psi^{-1})$ denotes the domain of definition for $\Psi^{-1}$.
\el

\section{The LDP for path-dependent multivalued  McKean-Vlasov SDEs}\label{sec:ldp}

In this section, we take ${\rm D}=\overline{\cD(A)}$, i.e. $\sC=C([-r_0,0], \overline{\cD(A)})$, introduce a general criterion for the LDP and then prove the LDP for path-dependent multivalued McKean-Vlasov SDEs. 

\subsection{A general criterion for the LDP}

In this subsection, we present a general criterion for the LDP. 

Let $(\mS,\rho)$ be a Polish space. For each $0<\e<1$, let $X^{\e}(\cdot)$ be a $\mS$-valued random variable on $(\Omega,\sF,\{\sF_t\}_{t\in[0,T]},\mP)$.

\bd
The function $I$ on $\mS$ is called a rate function if for each $M<\infty$, $\{\varsigma\in\mS: I(\varsigma)\leq M\}$ is a compact subset of $\mS$.
\ed

\bd
$\{X^{\e}(\cdot), \e\in(0,1)\}$ is said to satisfy the LDP with the speed $\e$ and the rate function $I$, if for any subset $O\in\sB(\mS)$,
$$
-\inf_{\varsigma\in {\rm Int}(O)}I(\varsigma)\leq\liminf_{\e\to0}\e\log\mP(X^{\e}(\cdot)\in {\rm Int}(O))\leq\limsup_{\e\to0}\e\log\mP(X^{\e}(\cdot)\in \bar{O})\leq -\inf_{\varsigma\in \bar{O}}I(\varsigma),
$$
where ${\rm Int}(O)$ and $\bar{O}$ denote the interior and the closure of $O$, respectively and they are taken in $\mS$.
\ed

\bd
$\{X^{\e}(\cdot), \e\in(0,1)\}$ is said to satisfy the Laplace principle with the speed $\e$ and the rate function $I$,  if for any real bounded continuous function $G$ on $\mathbb{S}$,
\ce
\lim\limits_{\varepsilon\rightarrow 0}\varepsilon \log \mE\left\{\exp\left[-\frac{G(X^{\e}(\cdot))}{\e}\right]\right\}=-\inf\limits_{\varsigma\in \mathbb{S}}\(G(\varsigma)+I(\varsigma)\).
\de
\ed

In the sequel, we take $\mS:=C([-r_0,T],\overline{\cD(A)})$ and $\rho(X^1, X^2):=\sup\limits_{t\in[-r_0,T]}|X^1(t)-X^2(t)|$ for $X^1, X^2\in C([-r_0,T],\overline{\cD(A)})$. Note that the LDP is equivalent to the Laplace principle (cf. \cite{de}). Therefore, in order to obtain the LDP for $\{X^{\e}(\cdot), \e\in(0,1)\}$, we prove the Laplace principle for $\{X^{\e}(\cdot), \e\in(0,1)\}$. Then we state the conditions under which the Laplace principle holds. 

For each $N\in \mN$, define 
$$
S^{N}:=\left\{h\in L^2([0,T],\mR^m):\int_{0}^{T}|h(s)|^2\dif s\leq N\right\}.
$$
Equipped with the weak topology, $S^N$ is a compact subset of $L^2([0,T], \mR^m)$. Let 
$$
\cS^N:=\left\{u:\Omega\times[0,T]\to\mR^m, u\ \text{is}\ \{\sF_t\}_{t\in[0,T]}\ \text{-predictable and }   u(\omega,\cdot)\in S^N\ \text{for }a.e. \omega\in\Omega\right\}.
$$

\bco\label{cond}
Let $\Gamma^{\e}:C([0,T], \mR^m)\to \mS$ be a family of measurable mappings. There exists a measurable mapping $\Gamma^0: C([0,T], \mR^m)\to \mS$ such that
\begin{enumerate}[$(i)$]
\item For each $N\in\mN$, any $\{h_\e\}\subset S^N$ and $h\in S^N$ satisfying $h_\e\to h$ as $\e\rightarrow0$, 
$$
\lim\limits_{\e\to0}\sup\limits_{t\in[-r_0,T]}\left|\Gamma^0\left(\int_{0}^{\cdot}h_\e(s)\dif s\right)(t)-\Gamma^{0}\left(\int_{0}^{\cdot}h(s)\dif s\right)(t)\right|=0.
$$
\item For each $N\in\mN$, any $\{u_\e\}\subset \cS^N$ and any $\delta>0$, 
$$
\lim_{\e\to0}\mP\left(\sup\limits_{t\in[-r_0,T]}\left|\Gamma^{\e}\left(W+\frac{1}{\sqrt{\e}}\int_{0}^{\cdot}u_{\e}(s)\dif s\right)(t)-\Gamma^{0}\left(\int_{0}^{\cdot}u_{\e}(s)\dif s\right)(t)\right|\geq \delta\right)=0.
$$
\end{enumerate}
\eco

Under the above assumptions, we have the following result (\cite[Theorem 3.2]{aw}).

\bt\label{ldpcon}
Set $X^{\e}(\cdot):=\Gamma^{\e}(W)$. Assume that $(i)$ and $(ii)$ in Condition \ref{cond} hold. Then the family $\{X^{\e}(\cdot), \e\in(0,1)\}$ satisfies the LDP in $\mS$ with the rate function $I$ given by
$$
I(g)=\frac{1}{2}\inf_{\substack{\{u\in L^2([0,T],\mR^m):\\ g=\Gamma^0(\int_{0}^{\cdot}u(s)\dif s)\}}}\left\{\int_{0}^{T}|u(s)|^2\dif s\right\},\quad g\in \mS,
$$
with the convention that $\inf\emptyset=+\infty$.
\et

\subsection{The LDP}

In this subsection, we prove the LDP for path-dependent multivalued McKean-Vlasov SDEs. 

We recall Eq.(\ref{eq0}), i.e.
\ce\left\{\begin{array}{ll}
	\dif X^{\e}(t)\in -A(X^{\e}(t))\dif t+b(X^{\e}_t,\sL_{X^{\e}_t})\dif t+\sqrt{\e}\s(X^{\e}_t,\sL_{X^{\e}_t})\dif W(t), t\in(0,T],\\
	X^\e(t)=\xi(t) \in \overline{\cD(A)}, t\in[-r_0,0].
\end{array}
\right.
\de

We make the following assumptions:
\begin{enumerate}[$({\bf H}_1)$] 
\item There exists a constant $L_1>0$ such that for any $\zeta\in\sC$, $\mu\in\sP_2^{\sC}$
$$
|b(\zeta,\mu)|^2+\|\s(\zeta,\mu)\|^2\leq L_1(1+\|\zeta\|^2_{\infty}+\mu(\|\cdot\|^2_{\infty})).
$$ 
\item $b(\cdot,\cdot)$ is continuous on $\sC\times \sP_{2}^{\sC}$, and there exists a constant $L_2>0$ such that for any $\zeta, \eta\in\sC$, $\mu,\nu\in\sP_2^{\sC}$
\ce
&&2\<b(\zeta,\mu)-b(\eta,\nu),\zeta(0)-\eta(0)\>\leq L_2\left(\kappa_1\(\|\zeta-\eta\|^2_\infty\)+\kappa_2\(\mW_2(\mu,\nu)^2\)\right),\\
&&\|\s(\zeta,\mu)-\s(\eta,\nu)\|^2\leq L_2\left(\kappa_1\(\|\zeta-\eta\|^2_\infty\)+\kappa_2\(\mW_2(\mu,\nu)^2\)\right),
\de
where $\kappa_i: [0, \infty)\to [0, \infty), i=1,2$ are strictly increasing, continuous concave functions and satisfy
\ce
\kappa_1(0)=\kappa_2(0)=0, \quad \int_{0+}\frac{1}{\kappa_1(v)+\kappa_2(v)}\dif v=\infty.
\de
\end{enumerate}

Under $({\bf H}_1)$ and $({\bf H}_2)$, by Theorem 4.1 in \cite{mq}, we know that Eq.(\ref{eq0}) has a unique strong solution denoted as $(X^{\e}_{\cdot},K^{\e}(\cdot))$ and there exists a functional $\Gamma^{\e}: C([0,T], \mR^m)\to C([-r_0,T],\overline{\cD(A)})$ such that $X^{\e}(\cdot)=\Gamma^{\e}(W)$. We first consider a controlled analogue of Eq.(\ref{eq0}) with the same initial value
\be\left\{\begin{array}{ll}
\dif X^{\e,u}(t)\in-A(X^{\e,u}(t))\dif t+b(X^{\e,u}_t,\sL_{X^{\e}_t})\dif t+\sigma(X^{\e,u}_t,\sL_{X^{\e}_t})u(t)\dif t\\
\qquad\qquad\qquad+\sqrt{\e}\sigma(X^{\e,u}_t,\sL_{X^{\e}_t})\dif W(t), \quad u\in\cS^N,  t\in(0,T],\\
X^{\e,u}(t)=\xi(t) \in \overline{\cD(A)}, t\in[-r_0,0].
\end{array}
\right.
\label{controlled}
\ee
By the Girsanov theorem, Eq.(\ref{controlled}) has a unique solution denoted by $(X^{\e,u}_\cdot,K^{\e,u}(\cdot))$. Moreover, we have $X^{\e,u}(\cdot)=\Gamma^{\e}(W+\frac{1}{\sqrt{\e}}\int_0^{\cdot}u(s)\dif s)$.

Let $(X^{0,u}_\cdot,K^{0,u}(\cdot))$ be the solution of the following deterministic equation
\be\left\{\begin{array}{ll}
\dif X^{0,u}(t)\in -A(X^{0,u}(t))\dif t+b(X^{0,u}_t,\delta_{X^{0}_t})\dif t+\sigma(X^{0,u}_t,\delta_{X^{0}_t})u(t)\dif t, u\in\cS^N, t\in(0,T],\\
X^{0,u}(t)=\xi(t) \in \overline{\cD(A)}, \quad t\in[-r_0,0],
\end{array}
\right.
\label{eqlim}
\ee
where $\delta_{X^{0}_t}$ is the Dirac measure at $X^{0}_t$ and $(X^0_\cdot,K^0(\cdot))$ is the solution of the following equation
\be\left\{\begin{array}{ll}
\dif X^0(t)\in -A(X^{0}(t))\dif t+b(X^0_t,\delta_{X^0_t})\dif t, \quad t\in(0,T],\\
X^{0}(t)=\xi(t) \in \overline{\cD(A)},\quad t\in[-r_0,0].
\end{array}
\right.
\label{eqct}
\ee
Define the measurable map $\Gamma^0:C([0,T],\mR^m)\to C([-r_0,T],\overline{\cD(A)})$ by $\Gamma^0(\int_{0}^{\cdot}u(s)\dif s)=X^{0,u}(\cdot)$.

Next, we justify $(i)$ and $(ii)$ in Condition \ref{cond} by $\Gamma^{\e}, \Gamma^0$.

\subsubsection{Verification of Condition \ref{cond} $(i)$}

First of all, we prepare two following moment estimates.

\bl
Suppose that $({\bf H}_1)$ and $({\bf H}_2)$ hold. Then Eq.(\ref{eqct}) has a unique solution $(X^0_\cdot,K^0(\cdot))$. Moreover, it holds that for any $v\in \mathrm{Int}(\cD(A))$
\be
\sup\limits_{r\in[-r_0,T]}|X^0(r)-v|^2\leq C_1e^{C_2T},
\label{x02moes}
\ee
where $C_1= 2\sup\limits_{r\in[-r_0,0]}|\xi(r)-v|^2+(2\g_2 +2\g_3 +L_1+4L_1|v|^2)T$, and $C_2=2\g_2+1+4L_1$.
\el
\begin{proof}
First of all, we follow the line of the proof in \cite[Theorem 4.1]{mq} and obtain that Eq.(\ref{eqct}) has a unique solution $(X^0_\cdot,K^0(\cdot))$.

Next, our aim is to show (\ref{x02moes}). By Lemma \ref{L2} and $({\bf H}_1)$, it holds that for any $v\in \mathrm{Int}(\cD(A))$ and $t\in[0,T]$,
\ce
|X^0(t)-v|^2&=&|\xi(0)-v|^2-2\int_{0}^{t}\<X^{0}(s)-v,\dif K^{0}(s)\>+2\int_{0}^{t}\<X^{0}(s)-v,b(X^0_s,\delta_{X^0_s})\>\dif s\\
&\leq&|\xi(0)-v|^2-2\g_1|K^{0}|_0^t+2\g_2\int_{0}^{t}|X^{0}(s)-v|\dif s+2\g_3 t+\int_{0}^{t}|X^{0}(s)-v|^2\dif s\\
&&+\int_{0}^{t}|b(X^0_s,\delta_{X^0_s})|^2\dif s\\
&\leq&|\xi(0)-v|^2+2\g_2\int_{0}^{t}(1+|X^{0}(s)-v|^2)\dif s+2\g_3 t+\int_{0}^{t}|X^{0}(s)-v|^2\dif s\\
&&+L_1\int_{0}^{t}(1+\|X^0_s\|^2_\infty+\delta_{X^0_s}(\|\cdot\|^2_\infty))\dif s.
\de
Note that $\delta_{X^0_s}(\|\cdot\|^2_\infty)=\|X^0_s\|^2_\infty$ and
\be
\|X^0_s\|^2_\infty=\sup\limits_{\t\in[-r_0,0]}|X^0(s+\t)|^2\leq\sup\limits_{r\in[-r_0,s]}|X^0(r)|^2\leq 2\sup\limits_{r\in[-r_0,s]}|X^0(r)-v|^2+2|v|^2.
\label{x0d0es}
\ee
Thus, we infer that
\ce
|X^0(t)-v|^2&\leq& |\xi(0)-v|^2+2\g_2 t+2\g_3 t+L_1t+(2\g_2+1)\int_{0}^{t}\sup\limits_{r\in[-r_0,s]}|X^{0}(r)-v|^2\dif s\\
&&+4L_1\int_0^t\sup\limits_{r\in[-r_0,s]}|X^0(r)-v|^2\dif s+4L_1|v|^2t,
\de
and
\ce
\sup\limits_{r\in[0,t]}|X^0(r)-v|^2&\leq& |\xi(0)-v|^2+2\g_2 t+2\g_3 t+L_1t+(2\g_2+1)\int_{0}^{t}\sup\limits_{r\in[-r_0,s]}|X^{0}(r)-v|^2\dif s\\
&&+4L_1\int_0^t\sup\limits_{r\in[-r_0,s]}|X^0(r)-v|^2\dif s+4L_1|v|^2t.
\de
Since 
\ce
\sup\limits_{r\in[-r_0,t]}|X^0(r)-v|^2\leq \sup\limits_{r\in[-r_0,0]}|X^0(r)-v|^2+\sup\limits_{r\in[0,t]}|X^0(r)-v|^2,
\de
and $X^0(r)=\xi(r), r\in[-r_0,0]$, it holds that
\ce
\sup\limits_{r\in[-r_0,t]}|X^0(r)-v|^2\leq \sup\limits_{r\in[-r_0,0]}|\xi(r)-v|^2+\sup\limits_{r\in[0,t]}|X^0(r)-v|^2.
\de
So, one can conclude that
\ce
\sup\limits_{r\in[-r_0,t]}|X^0(r)-v|^2&\leq& 2\sup\limits_{r\in[-r_0,0]}|\xi(r)-v|^2+(2\g_2 +2\g_3 +L_1+4L_1|v|^2)T\\
&&+(2\g_2+1+4L_1)\int_{0}^{t}\sup\limits_{r\in[-r_0,s]}|X^{0}(r)-v|^2\dif s.
\de
The Gronwall inequality yields the required estimate.
\end{proof}

\bl\label{ldplmlimbnd}
Assume that $({\bf H}_1)$ and $({\bf H}_2)$ hold. Then Eq.(\ref{eqlim}) has a unique solution $(X^{0,u}_\cdot,K^{0,u}(\cdot))$. Furthermore, it holds that
$$
\sup_{u\in \cS^N}\sup_{t\in[-r_0,T]}|X^{0,u}(t)|^2+\sup_{u\in \cS^N}|K^{0,u}|_0^T\leq C, \quad a.s.
$$
where the constant $C>0$ is independent of $\omega$.
\el
\begin{proof}
Firstly, based on the similar proof to that in \cite[Theorem 3.6]{mq}, one can obtain that Eq.(\ref{eqlim}) has a unique solution $(X^{0,u}_\cdot,K^{0,u}(\cdot))$. 

Next, by Lemma \ref{L2} and $({\bf H}_1)$, we have that for any $v\in \mathrm{Int}(\cD(A))$ and $t\in[0,T]$,
\ce
|X^{0,u}(t)-v|^2&=&|\xi(0)-v|^2-2\int_{0}^{t}\<X^{0,u}(s)-v,\dif K^{0,u}(s)\>\\
&&+2\int_{0}^{t}\<X^{0,u}(s)-v,b(X^{0,u}_s,\delta_{X^{0}_s})\>\dif s\\
&&+2\int_{0}^{t}\<X^{0,u}(s)-v,\s(X^{0,u}_s,\delta_{X^{0}_s})u(s)\>\dif s\\
&\leq&|\xi(0)-v|^2-2\gamma_1|K^{0,u}|_0^t+2\gamma_2\int_{0}^{t}|X^{0,u}(s)-v|\dif s+2\gamma_3t\\
&&+\int_{0}^{t}|X^{0,u}(s)-v|^2\dif s+\int_{0}^{t}|b(X^{0,u}_s,\delta_{X^{0}_s})|^2\dif s\\
&&+2\sup_{s\in[0,t]}|X^{0,u}(s)-v|\(\int_{0}^{t}\|\s(X^{0,u}_s,\delta_{X^{0}_s})\|^2\dif s\)^{1/2}\(\int_{0}^{t}|u(s)|^2\dif s\)^{1/2}\\
&\leq&|\xi(0)-v|^2-2\gamma_1|K^{0,u}|_0^t+(2\gamma_2+1)\int_{0}^{t}|X^{0,u}(s)-v|^2\dif s\\
&&+2\gamma_2t+2\gamma_3t+\int_{0}^{t}|b(X^{0,u}_s,\delta_{X^{0}_s})|^2\dif s\\
&&+\frac{1}{2}\sup_{s\in[0,t]}|X^{0,u}(s)-v|^2+2N\int_{0}^{t}\|\s(X^{0,u}_s,\delta_{X^{0}_s})\|^2\dif s, 
\de
where we use $u\in\cS^N$. Moreover, from $({\bf H}_1)$, it follows that
\ce
|X^{0,u}(t)-v|^2&\leq&|\xi(0)-v|^2-2\gamma_1|K^{0,u}|_0^t+(2\gamma_2+1)\int_{0}^{t}|X^{0,u}(s)-v|^2\dif s+2\gamma_2t+2\gamma_3t\\
&&+L_1(1+2N)\int_{0}^{t}\(1+\|X^{0,u}_s\|^2_{\infty}+\delta_{X^{0}_s}(\|\cdot\|^2_{\infty})\)\dif s\\
&&+\frac{1}{2}\sup_{s\in[0,t]}|X^{0,u}(s)-v|^2\\
&\leq&|\xi(0)-v|^2-2\gamma_1|K^{0,u}|_0^t+(2\gamma_2+1)\int_{0}^{t}|X^{0,u}(s)-v|^2\dif s+2\gamma_2t+2\gamma_3t\\
&&+L_1(1+2N)t+L_1(1+2N)\int_{0}^{t}\|X^{0,u}_s\|^2_{\infty}\dif s+L_1(1+2N)\int_{0}^{t}\|X^{0}_s\|^2_\infty\dif s\\
&&+\frac{1}{2}\sup_{s\in[0,t]}|X^{0,u}(s)-v|^2.
\de 
Note that
\ce
&&\|X^{0,u}_s\|^2_{\infty}\leq\sup_{r\in[-r_0,s]}|X^{0,u}(r)|^2\leq2\sup_{r\in[-r_0,s]}|X^{0,u}(r)-v|^2+2|v|^2,\\
&&\|X^{0}_s\|^2_{\infty}\leq\sup_{r\in[-r_0,s]}|X^{0}(r)|^2\leq2\sup_{r\in[-r_0,s]}|X^{0}(r)-v|^2+2|v|^2\leq 2C_1e^{C_2T}+2|v|^2,
\de
where we use (\ref{x02moes}) in the last inequality. Thus, we get that
\ce
|X^{0,u}(t)-v|^2+2\gamma_1|K^{0,u}|_0^t&\leq&|\xi(0)-v|^2+(2\gamma_2+1)\int_{0}^{t}|X^{0,u}(s)-v|^2\dif s+2\gamma_2t+2\gamma_3t\\
&&+L_1(1+2N)t+2L_1(1+2N)\int_{0}^{t}\sup_{r\in[-r_0,s]}|X^{0,u}(r)-v|^2\dif s\\
&&+2L_1(1+2N)|v|^2t+Ct+\frac{1}{2}\sup_{s\in[0,t]}|X^{0,u}(s)-v|^2,
\de
and furthermore,
\be
&&\sup_{r\in[0,t]}|X^{0,u}(r)-v|^2+4\gamma_1|K^{0,u}|_0^t\no\\
&\leq&2|\xi(0)-v|^2+2\(2\gamma_2+2\gamma_3+L_1(1+2N)+2L_1(1+2N)|v|^2+C\)T\no\\
&&+2(2\gamma_2+1+2L_1(1+2N))\int_{0}^{t}\sup_{r\in[-r_0,s]}|X^{0,u}(r)-v|^2\dif s.
\label{xesti}
\ee
Since $X^{0,u}(r)=\xi(r)$ for $r\in[-r_0,0]$, it holds that
\ce
\sup_{r\in[-r_0,t]}|X^{0,u}(r)-v|^2&\leq&\sup_{r\in[0,t]}|X^{0,u}(r)-v|^2+\sup_{r\in[-r_0,0]}|X^{0,u}(r)-v|^2\\
&=&\sup_{r\in[0,t]}|X^{0,u}(r)-v|^2+\sup_{r\in[-r_0,0]}|\xi(r)-v|^2.
\de
So, in terms of (\ref{xesti}), we conclude that
\ce
&&\sup_{r\in[-r_0,t]}|X^{0,u}(r)-v|^2\\
&\leq&3\sup_{r\in[-r_0,0]}|\xi(r)-v|^2+2\(2\gamma_2+2\gamma_3+L_1(1+2N)+2L_1(1+2N)|v|^2+C\)T\\
&&+2(2\gamma_2+1+2L_1(1+2N))\int_{0}^{t}\sup_{r\in[-r_0,s]}|X^{0,u}(r)-v|^2\dif s,
\de
which together with the Gronwall inequality yields that
\be
\sup_{r\in[-r_0,T]}|X^{0,u}(r)-v|^2\leq C'_1e^{C'_2T}, 
\label{x0hes}
\ee
where $C'_1=3\sup\limits_{r\in[-r_0,0]}|\xi(r)-v|^2+2\(2\gamma_2+2\gamma_3+L_1(1+2N)+2L_1(1+2N)|v|^2+C\)T$ and $C'_2=2(2\gamma_2+1+2L_1(1+2N))$.

Finally, by (\ref{x0hes}) it holds that
\ce
\sup_{t\in[-r_0,T]}|X^{0,u}(t)|^2\leq 2C'_1e^{C'_2T}+2|v|^2,
\de
and further
\ce
\sup_{u\in\cS^N}\sup_{t\in[-r_0,T]}|X^{0,u}(t)|^2\leq 2C'_1e^{C'_2T}+2|v|^2<\infty.
\de
Besides, from (\ref{xesti}) and (\ref{x0hes}), it follows that
\ce
\sup_{u\in\cS^N}|K^{0,u}|_0^T<\infty.
\de
The proof is complete. 
\end{proof}

The following limit result is also needed in the sequel.

\bl\label{hehconv}
Assume that $({\bf H}_1)$ and $({\bf H}_2)$ hold. For each $N\in\mN$, any $\{h_\e\}\subset S^N$ and $h\in S^N$, set
$$
g_\e(t):=\int_{0}^{t}\s(X^{0,h}_s,\delta_{X^0_s})(h_\e(s)-h(s))\dif s, \quad t\in[0,T].
$$
If $h_\e\to h$ in $S^N$ as $\e\rightarrow0$, then $g_\e\to0$ in $C([0,T], \mR^d)$. 
\el
\begin{proof}
First of all, we show that $\{g_\e, \e\in(0,1)\}$ is compact in $C([0,T], \mR^d)$. Based on the Arzel\`a-Ascoli theorem, we only need to prove that 
\begin{enumerate}[$(i)$]
\item $\{g_\e(t)\}$ is uniformly bounded, that is, $\sup\limits_{\e}\sup\limits_{t\in[0,T]}|g_\e(t)|<\infty$,
\item $\{g_\e(t)\}$ is equicontinuous in $C([0,T], \mR^d)$.
\end{enumerate}
Indeed, by $({\bf H}_1)$, it holds that
\ce
|g_\e(t)|&\leq&\int_{0}^{t}|\s(X^{0,h}_s,\delta_{X^0_s})(h_\e(s)-h(s)|\dif s\\
&\leq& \(\int_{0}^{t}\|\s(X^{0,h}_s,\delta_{X^0_s})\|^2\dif s\)^{1/2}\(\int_{0}^{t}|h_\e(s)-h(s)|^2\dif s\)^{1/2}\\
&\leq& \sqrt{L_1}\(\int_{0}^{t}(1+\|X^{0,h}_s\|_\infty^2+\delta_{X^0_s}(\|\cdot\|_\infty^2))\dif s\)^{1/2}\(2\int_{0}^{t}(|h_\e(s)|^2+|h(s)|^2)\dif s\)^{1/2}\\
&\leq&2\sqrt{L_1N T}\(1+\sup_{t\in[-r_0,T]}(|X^{0,h}(t)|^2+|X^0(t)|^2)\)^{1/2},
\de
which together with (\ref{x02moes}) and Lemma \ref{ldplmlimbnd} yields that
\ce
\sup\limits_{\e}\sup_{t\in[0,T]}|g_\e(t)|<\infty.
\de
That is, $(i)$ is right.

Further, by the similar deduction to the above, we have that for any $\d>0$
\ce
|g_\e(t+\d)-g_\e(t)|&\leq& \int_{t}^{t+\d}|\s(X^{0,h}_s,\delta_{X^0_s})(h_\e(s)-h(s))|\dif s\\
&\leq& 2\sqrt{L_1N\d}\(1+\sup_{t\in[-r_0,T]}(|X^{0,h}(t)|^2+|X^0(t)|^2)\)^{1/2},
\de
and
\ce
\lim_{\d\to 0}\sup_{\e}|g_\e(t+\d)-g_\e(t)|=0, \quad t\in[0,T],
\de
which implies $(ii)$. 

Next, we notice that for any $\varpi\in\mR^d$,
\ce
\int_0^T|\s^*(X^{0,h}_s,\delta_{X^0_s})\varpi|^2\dif s&\leq& \int_0^T\|\s^*(X^{0,h}_s,\delta_{X^0_s})\|^2|\varpi|^2\dif s\\
&\leq& L_1|\varpi|^2 \int_0^T(1+\|X^{0,h}_s\|_\infty^2+\delta_{X^0_s}(\|\cdot\|_\infty^2))\dif s\\
&\leq& L_1|\varpi|^2T\(1+\sup_{t\in[-r_0,T]}(|X^{0,h}(t)|^2+|X^0(t)|^2)\)<\infty,
\de
where we use $({\bf H}_1)$, (\ref{x02moes}) and Lemma \ref{ldplmlimbnd}. So, when $h_\e$ converges weakly to $h$ in $L^2([0,T], \mR^m)$ as $\e\to0$,
\ce
\int_{0}^{t}\<h_\e(s),\s^*(X^{0,h}_s,\delta_{X^0_s})\varpi\>\dif s\to \int_{0}^{t}\<h(s),\s^*(X^{0,h}_s,\delta_{X^0_s})\varpi\>\dif s, 
\de
and
\ce
\int_{0}^{t}\s(X^{0,h}_s,\delta_{X^0_s})h_\e(s)\dif s\to \int_{0}^{t}\s(X^{0,h}_s,\delta_{X^0_s})h(s)\dif s,
\de
which together with the compactness of $\{g_\e, \e\in(0,1)\}$ in $C([0,T], \mR^d)$ implies that
\ce
\lim_{\e\to0}\sup_{t\in[0,T]}|g_\e(t)|=0.
\de
Thus, the proof is complete.
\end{proof}

Now it is the position to verify Condition \ref{cond} $(i)$.

\bp\label{ldpco1}
Assume $({\bf H}_1)$ and $({\bf H}_2)$ hold. For each $N\in\mN$, any $\{h_\e\}\subset S^N$ and $h\in S^N$ satisfying $h_\e\to h$ as $\e\rightarrow0$, it holds that
$$
\lim\limits_{\e\to0}\sup\limits_{t\in[-r_0,T]}\left|\Gamma^0\left(\int_{0}^{\cdot}h_\e(s)\dif s\right)(t)-\Gamma^{0}\left(\int_{0}^{\cdot}h(s)\dif s\right)(t)\right|=0.
$$
\ep
\begin{proof}
By the definition of $\Gamma^0$, $X^{0,h}(\cdot)=\Gamma^0(\int_{0}^{\cdot}h(s)\dif s)(\cdot)$ is the solution of Eq.(\ref{eqlim}) and  $X^{0,h_\e}(\cdot)=\Gamma^0(\int_{0}^{\cdot}h_\e(s)\dif s)(\cdot)$ is the solution of Eq.(\ref{eqlim}) with $h$ replaced by $h_\e$. Therefore, in order to obtain the required limit, we only need to estimate $\sup\limits_{t\in[-r_0,T]}|X^{0,h_\e}(t)-X^{0,h}(t)|$.

Set $Z^0(t):=X^{0,h_\e}(t)-X^{0,h}(t)$. By Lemma \ref{L1}, it holds that
\ce
|Z^0(t)|^2&=&-2\int_{0}^{t}\<Z^0(s),\dif K^{0,h_\e}(s)-\dif K^{0,h}(s)\>\\
&&+2\int_{0}^{t}\<Z^0(s),b(X^{0,h_\e}_s,\delta_{X^0_s})-b(X^{0,h}_s,\delta_{X^0_s})\>\dif s\\
&&+2\int_{0}^{t}\<Z^0(s),\s(X^{0,h_\e}_s,\delta_{X^0_s})h_\e(s)-\s(X^{0,h}_s,\delta_{X^0_s})h(s)\>\dif s\\
&\leq&2\int_{0}^{t}\<Z^0(s),b(X^{0,h_\e}_s,\delta_{X^0_s})-b(X^{0,h}_s,\delta_{X^0_s})\>\dif s\\
&&+2\int_{0}^{t}\<Z^0(s),\(\s(X^{0,h_\e}_s,\delta_{X^0_s})-\s(X^{0,h}_s,\delta_{X^0_s})\)h_\e(s)\>\dif s\\
&&+2\int_{0}^{t}\<Z^0(s),\s(X^{0,h}_s,\delta_{X^0_s})(h_\e(s)-h(s))\>\dif s\\
&=:&J_1(t)+J_2(t)+J_3(t).
\de
For $J_1(t)$, according to $({\bf H}_2)$, we have that
\ce
J_1(t)\leq L_2\int_{0}^{t}\kappa_1(\|X^{0,h_\e}_s-X^{0,h}_s\|^2_{\infty})\dif s. 
\de
For $J_2(t)$, from the Young inequality, the H$\ddot{\rm o}$lder inequality and $h_\e\in S^N$, it follows that
\ce
J_2(t)&\leq&2\int_{0}^{t}|Z^0(s)|\|\s(X^{0,h_\e}_s,\delta_{X^0_s})-\s(X^{0,h}_s,\delta_{X^0_s})\||h_\e(s)|\dif s\\
&\leq& 2\sup_{s\in[0,t]}|Z^0(s)|\int_{0}^{t}\|\s(X^{0,h_\e}_s,\delta_{X^0_s})-\s(X^{0,h}_s,\delta_{X^0_s})\||h_\e(s)|\dif s\\
&\leq&2\(\int_{0}^{t}\|\s(X^{0,h_\e}_s,\delta_{X^0_s})-\s(X^{0,h}_s,\delta_{X^0_s})\|^2\dif s\)\(\int_{0}^{t}|h_\e(s)|^2\dif s\)\\
&&+\frac{1}{2}\sup_{s\in[0,t]}|Z^0(s)|^2\\
&\leq&2NL_2\int_{0}^{t}\kappa_1(\|X^{0,h_\e}_s-X^{0,h}_s\|_\infty^2)\dif s+\frac{1}{2}\sup_{s\in[0,t]}|Z^0(s)|^2.
\de
So, the above deduction implies that
\ce
|Z^0(t)|^2\leq(L_2+2NL_2)\int_{0}^{t}\kappa_1(\|X^{0,h_\e}_s-X^{0,h}_s\|^2_{\infty})\dif s+\frac{1}{2}\sup_{s\in[0,t]}|Z^0(s)|^2+|J_3(t)|,
\de
and
\ce
\sup\limits_{s\in[0,t]}|Z^0(s)|^2\leq 2(L_2+2NL_2)\int_{0}^{t}\kappa_1(\|X^{0,h_\e}_s-X^{0,h}_s\|^2_{\infty})\dif s+2\sup\limits_{s\in[0,t]}|J_3(s)|.
\de
Note that
\ce
\|X^{0,h_\e}_s-X^{0,h}_s\|^2_{\infty}=\sup\limits_{\t\in[-r_0,0]}|X^{0,h_\e}(s+\t)-X^{0,h}(s+\t)|^2\leq \sup\limits_{r\in[-r_0,s]}|X^{0,h_\e}(r)-X^{0,h}(r)|^2,
\de
and $X^{0,h_\e}(r)=X^{0,h}(r)=\xi(r)$ for $r\in[-r_0,0]$. Thus, it holds that
\ce
\|X^{0,h_\e}_s-X^{0,h}_s\|^2_{\infty}\leq \sup\limits_{r\in[0,s]}|X^{0,h_\e}(r)-X^{0,h}(r)|^2.
\de

Collecting the above deduction, we conclude that
\ce
\sup\limits_{s\in[0,t]}|Z^0(s)|^2&\leq&2(L_2+2NL_2)\int_{0}^{t}\kappa_1(\sup_{r\in[0,s]}|Z^0(r)|^2)\dif s+2\sup\limits_{s\in[0,T]}|J_3(s)|.
\de
Moreover, by Lemma \ref{bihari}, it holds that
\be
\sup_{s\in[0,T]}|Z^0(s)|^2\leq \Psi_1^{-1}\(\Psi_1(2\sup_{s\in[0,T]}|J_3(s)|)+2(L_2+2NL_2)T\),
\label{Zbi}
\ee
where $\Psi_1(r)=\int_\vartheta^r\frac{1}{\kappa_1(v)}\dif v$ for $\vartheta>0$ and $r\geq 0$. We claim:
\be
\lim\limits_{\e\to 0}\sup_{s\in[0,T]}|J_3(s)|=0.
\label{j30}
\ee
Taking the limit on two sides of (\ref{Zbi}) as $\e$ tends to $0$, by the property of $\kappa_1$ we obtain the required limit.

Next, we prove (\ref{j30}). By the Taylor formula, it holds that
\ce
\<Z^0(t),g_\e(t)\>&=&\int_{0}^{t}\<Z^0(t),\s(X^{0,h}_s,\delta_{X^0_s})[h_\e(s)-h(s)]\>\dif s-\int_{0}^{t}\<g_\e(s),\dif (K^{0,u_\e}(s)-K^{0,u}(s))\>\\
&&+\int_{0}^{t}\<g_\e(s),b(X^{0,h_\e}_s,\delta_{X^0_s})-b(X^{0,h}_s,\delta_{X^0_s})\>\dif s\\
&&+\int_{0}^{t}\<g_\e(s),\(\s(X^{0,h_\e}_s,\delta_{X^0_s})h_\e(s)-\s(X^{0,h}_s,\delta_{X^0_s})h(s)\)\>\dif s,
\de
which implies that
\ce
J_3(t)&=&2\<Z^0(t),g_\e(t)\>+2\int_{0}^{t}\<g_\e(s),\dif (K^{0,u_\e}(s)-K^{0,u}(s))\>\\
&&-2\int_{0}^{t}\<g_\e(s),b(X^{0,h_\e}_s,\delta_{X^0_s})-b(X^{0,h}_s,\delta_{X^0_s})\>\dif s\\
&&-2\int_{0}^{t}\<g_\e(s),\(\s(X^{0,h_\e}_s,\delta_{X^0_s})h_\e(s)-\s(X^{0,h}_s,\delta_{X^0_s})h(s)\)\>\dif s\\
&=:&J_{31}(t)+J_{32}(t)+J_{33}(t)+J_{34}(t).
\de
For $J_{31}(t)$, we notice that
\ce
\sup_{t\in[0,T]}|J_{31}(t)|\leq 2\sup_{t\in[0,T]}|g_\e(t)|\sup_{t\in[0,T]}|Z^0(t)|.
\de
So, by Lemma \ref{ldplmlimbnd} and \ref{hehconv}, it holds that
\ce
\lim_{\e\to0}\sup_{t\in[0,T]}|J_{31}(t)|=0.
\de
For $J_{32}(t)$, some simple computation implies that
\ce
\sup_{t\in[0,T]}|J_{32}(t)|&\leq&2\sup_{t\in[0,T]}|g_\e(t)|\(|K^{0,h_\e}|_0^T+|K^{0,h}|_0^T\),
\de
which together with Lemma \ref{ldplmlimbnd} and \ref{hehconv} yields that
\ce
\lim_{\e\to0}\sup_{t\in[0,T]}|J_{32}(t)|=0.
\de
For $J_{33}(t)$, by $({\bf H}_1)$ we get that
\ce
&&\sup_{t\in[0,T]}|J_{33}(t)|\\
&\leq&2\sup_{t\in[0,T]}|g_\e(t)|\int_{0}^{T}\(|b(X^{0,h_\e}_s,\delta_{X^0_s})|+|b(X^{0,h}_s,\delta_{X^0_s})|\)\dif s\no\\
&\leq&2\sup_{t\in[0,T]}|g_\e(t)|T^{1/2}\(\int_{0}^{T}2\(|b(X^{0,h_\e}_s,\delta_{X^0_s})|^2+|b(X^{0,h}_s,\delta_{X^0_s})|^2\)\dif s\)^{1/2}\no\\
&\leq&2\sup_{t\in[0,T]}|g_\e(t)|T^{1/2}\(\int_{0}^{T}2L_1\(2+\|X^{0,h_\e}_s\|^2_{\infty}+\|X^{0,h}_s\|^2_{\infty}+2\delta_{X^0_s}(\|\cdot\|^2_{\infty})\)\dif s\)^{1/2}.
\de
Then Lemma \ref{ldplmlimbnd} and \ref{hehconv} assure that
\ce
\lim_{\e\to0}\sup_{t\in[0,T]}|J_{33}(t)|=0.
\de
For $J_{34}(t)$, by some calculation, it holds that
\ce
&&\sup_{t\in[0,T]}|J_{34}(t)|\\
&\leq&2\sup_{t\in[0,T]}|g_\e(t)|\int_{0}^{T}\(\|\s(X^{0,h_\e}_s,\delta_{X^0_s})\||h_\e(s)|+\|\s(X^{0,h}_s,\delta_{X^0_s})\||h(s)|\)\dif s\\
&\leq&2\sup_{t\in[0,T]}|g_\e(t)|\[\(\int_{0}^{T}\|\s(X^{0,h_\e}_s,\delta_{X^0_s})\|^2\dif s\)^{1/2}\(\int_{0}^{T}|h_\e(s)|^2\dif s\)^{1/2}\\
&&+\(\int_{0}^{T}\|\s(X^{0,h}_s,\delta_{X^0_s})\|^2\dif s\)^{1/2}\(\int_{0}^{T}|h(s)|^2\dif s\)^{1/2}\]\\
&\leq&2\sup_{t\in[0,T]}|g_\e(t)|\sqrt{N}\[\(\int_{0}^{T}\|\s(X^{0,h_\e}_s,\delta_{X^0_s})\|^2\dif s\)^{1/2}+\(\int_{0}^{T}\|\s(X^{0,h}_s,\delta_{X^0_s})\|^2\dif s\)^{1/2}\].
\de
By $({\bf H}_1)$, (\ref{x02moes}) and Lemma \ref{ldplmlimbnd} and \ref{hehconv}, we obtain that
\ce
\lim_{\e\to0}\sup_{t\in[0,T]}|J_{34}(t)|=0.
\de
Finally, the above estimates give (\ref{j30}). The proof is complete.
\end{proof}

\subsubsection{Verification of Condition \ref{cond} $(ii)$}

First of all, we prepare the following lemma.

\bl\label{lmcenterdis}
Under $({\bf H}_1)$ and $({\bf H}_2)$, it holds that
$$
\lim_{\e\to0}\mE\(\sup_{t\in[-r_0,T]}|X^{\e}(t)-X^0(t)|^2\)=0.
$$
\el
\begin{proof}
By the It\^o formula and Lemma \ref{L1}, it holds that
\ce
|X^{\e}(t)-X^0(t)|^2&=&-2\int_{0}^{t}\<X^{\e}(s)-X^0(s),\dif \(K^{\e}(s)-K^0(s)\)\>\\
&&+2\int_{0}^{t}\<X^{\e}(s)-X^0(s),b(X^{\e}_s,\sL_{X^{\e}_s})-b(X^0_s,\delta_{X^0_s})\>\dif s\\
&&+2\sqrt{\e}\int_{0}^{t}\<X^{\e}(s)-X^0(s),\s(X^{\e}_s,\sL_{X^{\e}_s})\dif W(s)\>\\
&&+\e\int_{0}^{t}\|\s(X^{\e}_s,\sL_{X^{\e}_s})\|^2\dif s\\
&\leq&2\int_{0}^{t}\<X^{\e}(s)-X^0(s),b(X^{\e}_s,\sL_{X^{\e}_s})-b(X^0_s,\delta_{X^0_s})\>\dif s\\
&&+2\sqrt{\e}\int_{0}^{t}\<X^{\e}(s)-X^0(s),\s(X^{\e}_s,\sL_{X^{\e}_s})\dif W(s)\>\\
&&+\e\int_{0}^{t}\|\s(X^{\e}_s,\sL_{X^{\e}_s})\|^2\dif s\\
&=:&I_1(t)+I_2(t)+I_3(t).
\de
For $I_1(t)$, it follows from $({\bf H}_2)$ that
\ce
I_1(t)\leq L_2\int_{0}^{t}\(\kappa_1(\|X^{\e}_s-X^0_s\|^2_\infty)+\kappa_2(\mW_2(\sL_{X^{\e}_s},\delta_{X^0_s})^2)\)\dif s.
\de
In view of $I_3(t)$, by $({\bf H}_1)$ and $({\bf H}_2)$ we have that
\be
I_3(t)&\leq&2\e\int_{0}^{t}\|\s(X^{\e}_s,\sL_{X^{\e}_s})-\s(X^0_s,\delta_{X^0_s})\|^2\dif s+2\e\int_{0}^{t}\|\s(X^0_s,\delta_{X^0_s})\|^2\dif s\no\\
&\leq&2L_2\int_{0}^{t}\(\kappa_1(\|X^\e_s-X^0_s\|^2_{\infty})+\kappa_2(\mW_2(\sL_{X^{\e}_s},\delta_{X^0_s})^2)\)\dif s\no\\
&&+2\e L_1\int_{0}^{t}\(1+\|X^0_s\|^2_{\infty}+\delta_{X^0_s}(\|\cdot\|^2_\infty)\)\dif s,
\label{i3esti}
\ee
where we use $0<\e<1$. Collecting the above deduction, we conclude that
\ce
|X^{\e}(t)-X^0(t)|^2&\leq&3L_2\int_{0}^{t}\(\kappa_1(\|X^{\e}_s-X^0_s\|^2_\infty)+\kappa_2(\mW_2(\sL_{X^{\e}_s},\delta_{X^0_s})^2)\)\dif s\\
&&+2\e L_1\int_{0}^{t}\(1+\|X^0_s\|^2_{\infty}+\delta_{X^0_s}(\|\cdot\|^2_\infty)\)\dif s+|I_2(t)|,
\de
and
\ce
\sup\limits_{s\in[0,t]}|X^{\e}(s)-X^0(s)|^2&\leq&3L_2\int_{0}^{t}\(\kappa_1(\|X^{\e}_s-X^0_s\|^2_\infty)+\kappa_2(\mW_2(\sL_{X^{\e}_s},\delta_{X^0_s})^2)\)\dif s\\
&&+2\e L_1\int_{0}^{t}\(1+\|X^0_s\|^2_{\infty}+\delta_{X^0_s}(\|\cdot\|^2_\infty)\)\dif s+\sup\limits_{s\in[0,t]}|I_2(s)|.
\de
Besides, by the Burkholder-Davis-Gundy inequality and the similar deduction to that for (\ref{i3esti}), it holds that
\ce
&&\mE\(\sup_{s\in[0,t]}|I_2(s)|\)\\
&\leq& 2\sqrt{\e}C\mE\left(\int_{0}^{t}|X^{\e}(s)-X^0(s)|^2\|\s(X^{\e}_s,\sL_{X^{\e}_s})\|^2 \dif s\right)^{1/2}\\
&\leq& 2\sqrt{\e}C\mE\sup_{s\in[0,t]}|X^{\e}(s)-X^0(s)|\left(\int_{0}^{t}\|\s(X^{\e}_s,\sL_{X^{\e}_s})\|^2\dif s\right)^{1/2}\\
&\leq& \frac{1}{2}\mE\sup_{s\in[0,t]}|X^{\e}(s)-X^0(s)|^2+C\e\mE\int_{0}^{t}\|\s(X^{\e}_s,\sL_{X^{\e}_s})\|^2\dif s\\
&\leq&\frac{1}{2}\mE\sup_{s\in[0,t]}|X^{\e}(s)-X^0(s)|^2+2L_2C\mE\int_{0}^{t}\(\kappa_1(\|X^\e_s-X^0_s\|^2_{\infty})+\kappa_2(\mW_2(\sL_{X^{\e}_s},\delta_{X^0_s})^2)\)\dif s\\
&&+2\e L_1C\int_{0}^{t}\(1+\|X^0_s\|^2_{\infty}+\delta_{X^0_s}(\|\cdot\|^2_\infty)\)\dif s.
\de
So, we infer that
\ce
\mE\sup\limits_{s\in[0,t]}|X^{\e}(s)-X^0(s)|^2&\leq&2(3L_2+2L_2C)\mE\int_{0}^{t}\(\kappa_1(\|X^{\e}_s-X^0_s\|^2_\infty)+\kappa_2(\mW_2(\sL_{X^{\e}_s},\delta_{X^0_s})^2)\)\dif s\\
&&+4\e L_1(C+1)\int_{0}^{T}\(1+\|X^0_s\|^2_{\infty}+\delta_{X^0_s}(\|\cdot\|^2_\infty)\)\dif s.
\de

Note that $X^\e(r)=X^0(r)=\xi(r)$ for any $r\in[-r_0,0]$. Thus, it holds that
\be
\mW_2(\sL_{X^{\e}_s},\delta_{X^0_s})^2&\leq&\mE\|X^{\e}_s-X^0_s\|^2_{\infty}=\mE\sup\limits_{\t\in[-r_0,0]}|X^{\e}(s+\t)-X^0(s+\t)|^2\no\\
&\leq&\mE\sup\limits_{r\in[0,s]}|X^{\e}(r)-X^0(r)|^2.
\label{xex0distes}
\ee
Besides, we know that
\ce
\delta_{X^0_s}(\|\cdot\|^2_\infty)=\|X^0_s\|^2_{\infty}=\sup\limits_{\t\in[-r_0,0]}|X^{0}(s+\t)|^2\leq \sup\limits_{r\in[-r_0,T]}|X^{0}(r)|^2.
\de
Therefore, the concavity of $\kappa_1$ and (\ref{x02moes}) imply that
\ce
&&\mE\sup\limits_{s\in[0,t]}|X^{\e}(s)-X^0(s)|^2\\
&\leq&2(3L_2+2L_2C)\int_{0}^{t}\(\kappa_1(\mE\sup\limits_{r\in[0,s]}|X^{\e}(r)-X^0(r)|^2)+\kappa_2(\mE\sup\limits_{r\in[0,s]}|X^{\e}(r)-X^0(r)|^2)\)\dif s\\
&&+\e C.
\de
By Lemma \ref{bihari}, we get that
\ce
\mE\sup\limits_{s\in[0,T]}|X^{\e}(s)-X^0(s)|^2&\leq& \Psi_2^{-1}(\Psi_2(\e C)+2(3L_2+2L_2C)T),
\de
where $\Psi_2(r)=\int_{\vartheta}^{r}\frac{1}{\kappa_1(v)+\kappa_2(v)}\dif v$ for $\vartheta>0, r\geq 0$. Lastly, taking the limit on two sides of the above inequality as $\e\to0$, by Lemma \ref{2k} we establish the required limit. 
\end{proof}

At present, we are ready to justify Condition \ref{cond} $(ii)$. 

\bp\label{ldpco2}
Suppose that $({\bf H}_1)$ and $({\bf H}_2)$ hold. Then for each $N\in\mN$, any $\{u_\e\}\subset \cS^N$ and any $\delta>0$, 
$$
\lim_{\e\to0}\mP\left(\sup\limits_{t\in[-r_0,T]}\left|\Gamma^{\e}\left(W+\frac{1}{\sqrt{\e}}\int_{0}^{\cdot}u_{\e}(s)\dif s\right)(t)-\Gamma^{0}\left(\int_{0}^{\cdot}u_{\e}(s)\dif s\right)(t)\right|\geq \delta\right)=0.
$$
\ep
\begin{proof}
We divide the proof into two steps. In the first step, we establish a key limiting result; in the second step, we demonstrate that the desired limit holds.

{\bf Step 1.} We prove that
$$
\lim_{\e\to0}\mE\left(\sup\limits_{t\in[-r_0,T]}\left|\Gamma^{\e}\left(W+\frac{1}{\sqrt{\e}}\int_{0}^{\cdot}u_{\e}(s)\dif s\right)(t)-\Gamma^{0}\left(\int_{0}^{\cdot}u_{\e}(s)\dif s\right)(t)\right|^2\right)=0.
$$

By the definition of $\Gamma^{\e}$ and $\Gamma^0$, we have that
\ce
X^{\e,u_{\e}}(\cdot)=\Gamma^{\e}\left(W+\frac{1}{\sqrt{\e}}\int_{0}^{\cdot}u_{\e}(s)\dif s\right),\quad X^{0,u_{\e}}(\cdot)=\Gamma^{0}\left(\int_{0}^{\cdot}u_{\e}(s)\dif s\right).
\de
Thus, it is sufficient to prove that
$$
\lim_{\e\to0}\mE\(\sup_{t\in[-r_0,T]}|X^{\e,u_{\e}}(t)-X^{0,u_{\e}}(t)|^2\)=0.
$$
By the It\^o formula, it holds that for any $t\in[0,T]$
\ce
&&|X^{\e,u_{\e}}(t)-X^{0,u_{\e}}(t)|^2\\
&=&-2\int_{0}^{t}\<X^{\e,u_{\e}}(s)-X^{0,u_{\e}}(s),\dif \(K^{\e,u_{\e}}(s)-K^{0,u_{\e}}(s)\)\>\\
&&+2\int_{0}^{t}\<X^{\e,u_{\e}}(s)-X^{0,u_{\e}}(s),b(X^{\e,u_{\e}}_s,\sL_{X^{\e}_s})-b(X^{0,u_{\e}}_s,\delta_{X^0_s})\>\dif s\\
&&+2\int_{0}^{t}\<X^{\e,u_{\e}}(s)-X^{0,u_{\e}}(s),\(\s(X^{\e,u_{\e}}_s,\sL_{X^{\e}_s})-\s(X^{0,u_{\e}}_s,\delta_{X^0_s})\)u_{\e}(s)\>\dif s\\
&&+2\sqrt{\e}\int_{0}^{t}\<X^{\e,u_{\e}}(s)-X^{0,u_{\e}}(s),\s(X^{\e,u_{\e}}_s,\sL_{X^{\e}_s})\dif W(s)\>\\
&&+\e\int_{0}^{t}\|\s(X^{\e,u_{\e}}_s,\sL_{X^{\e}_s})\|^2\dif s\\
&=:&J_1(t)+J_2(t)+J_3(t)+J_4(t)+J_5(t).
\de
According to Lemma \ref{L1}, it holds that
$$
J_1(t)\leq0.
$$
For $J_2(t)$, by $({\bf H}_2)$, we get that
\ce
J_2(t)\leq L_2\int_{0}^{t}\(\kappa_1(\|X^{\e,u_{\e}}_s-X^{0,u_{\e}}_s\|^2_{\infty})+\kappa_2(\mW_2(\sL_{X^{\e}_s},\delta_{X^0_s})^2)\)\dif s.
\de
For $J_3(t)$, the H$\ddot{\rm o}$lder inequality, $({\bf H}_2)$ and $u_{\e}\in\cS^N$ imply that
\ce
J_3(t)&\leq&2\int_{0}^{t}|X^{\e,u_{\e}}(s)-X^{0,u_{\e}}(s)|\|\s(X^{\e,u_{\e}}_s,\sL_{X^{\e}_s})-\s(X^{0,u_{\e}}_s,\delta_{X^0_s})\||u_{\e}(s)|\dif s\no\\
&\leq&2\sup_{s\in[0,t]}|X^{\e,u_{\e}}(s)-X^{0,u_{\e}}(s)|\int_{0}^{t}\|\s(X^{\e,u_{\e}}_s,\sL_{X^{\e}_s})-\s(X^{0,u_{\e}}_s,\delta_{X^0_s})\||u_{\e}(s)|\dif s\no\\
&\leq& \frac{1}{4}\sup_{s\in[0,t]}|X^{\e,u_{\e}}(s)-X^{0,u_{\e}}(s)|^2\no\\
&&+ 4\(\int_{0}^{t}\|\s(X^{\e,u_{\e}}_s,\sL_{X^{\e}_s})-\s(X^{0,u_{\e}}_s,\delta_{X^0_s})\|^2\dif s\)\(\int_{0}^{t}|u_{\e}(s)|^2\dif s\)\no\\
&\leq& \frac{1}{4}\sup_{s\in[0,t]}|X^{\e,u_{\e}}(s)-X^{0,u_{\e}}(s)|^2\no\\
&&+ 4NL_2\int_{0}^{t}\(\kappa_1(\|X^{\e,u_{\e}}_s-X^{0,u_{\e}}_s\|^2_{\infty})+\kappa_2(\mW_2(\sL_{X^{\e}_s},\delta_{X^0_s})^2)\)\dif s.
\de
As for $J_5(t)$, by $({\bf H}_1)$ and $({\bf H}_2)$, it holds that
\be
J_5(t)&\leq&2\int_{0}^{t}\|\s(X^{\e,u_{\e}}_s,\sL_{X^{\e}_s})-\s(X^{0,u_{\e}}_s,\delta_{X^0_s})\|^2\dif s+2\e\int_{0}^{t}\|\s(X^{0,u_{\e}}_s,\delta_{X^0_s})\|^2\dif s\no\\
&\leq&2L_2\int_{0}^{t}\(\kappa_1(\|X^{\e,u_{\e}}_s-X^{0,u_\e}_s\|^2_{\infty})+\kappa_2(\mW_2(\sL_{X^{\e}_s},\delta_{X^0_s})^2)\)\dif s\no\\
&&+2L_1\e\int_{0}^{t}\(1+\|X^{0,u_{\e}}_s\|^2_{\infty}+\delta_{X^0_s}(\|\cdot\|^2_{\infty})\)\dif s,
\label{j5estima}
\ee
where we make use of $0<\e<1$. Combining the above deduction, we conclude that
\ce
&&|X^{\e,u_{\e}}(t)-X^{0,u_{\e}}(t)|^2\\
&\leq&\frac{1}{4}\sup_{s\in[0,t]}|X^{\e,u_{\e}}(s)-X^{0,u_{\e}}(s)|^2+2L_1\e\int_{0}^{t}\(1+\|X^{0,u_{\e}}_s\|^2_{\infty}+\delta_{X^0_s}(\|\cdot\|^2_{\infty})\)\dif s\\
&&+(3L_2+4NL_2)\int_{0}^{t}\(\kappa_1(\|X^{\e,u_{\e}}_s-X^{0,u_\e}_s\|^2_{\infty})+\kappa_2(\mW_2(\sL_{X^{\e}_s},\delta_{X^0_s})^2)\)\dif s\\
&&+|J_4(t)|,
\de
and
\ce
&&\sup\limits_{s\in[0,t]}|X^{\e,u_{\e}}(s)-X^{0,u_{\e}}(s)|^2\\
&\leq&\frac{1}{4}\sup_{s\in[0,t]}|X^{\e,u_{\e}}(s)-X^{0,u_{\e}}(s)|^2+2L_1\e\int_{0}^{t}\(1+\|X^{0,u_{\e}}_s\|^2_{\infty}+\delta_{X^0_s}(\|\cdot\|^2_{\infty})\)\dif s\\
&&+(3L_2+4NL_2)\int_{0}^{t}\(\kappa_1(\|X^{\e,u_{\e}}_s-X^{0,u_\e}_s\|^2_{\infty})+\kappa_2(\mW_2(\sL_{X^{\e}_s},\delta_{X^0_s})^2)\)\dif s\\
&&+\sup\limits_{s\in[0,t]}|J_4(s)|.
\de

Next, by the Burkholder-Davis-Gundy inequality and the same deduction to that for (\ref{j5estima}), it holds that
\ce
&&\mE\(\sup_{s\in[0,t]}|J_4(s)|\)\no\\
&\leq&2\sqrt{\e}C\mE\(\int_{0}^{t}|X^{\e,u_{\e}}(s)-X^{0,u_{\e}}(s)|^2\|\s(X^{\e,u_{\e}}_s,\sL_{X^{\e}_s})\|^2\dif s\)^{1/2}\no\\
&\leq& 2\sqrt{\e}C\mE\sup_{s\in[0,t]}|X^{\e,u_{\e}}(s)-X^{0,u_{\e}}(s)|\(\int_{0}^{t}\|\s(X^{\e,u_{\e}}_s,\sL_{X^{\e}_s})\|^2\dif s\)^{1/2}\no\\
&\leq& \frac{1}{4}\mE\sup_{s\in[0,t]}|X^{\e,u_{\e}}(s)-X^{0,u_{\e}}(s)|^2+\e C\mE\int_{0}^{t}\|\s(X^{\e,u_{\e}}_s,\sL_{X^{\e}_s})\|^2\dif s\no\\
&\leq& \frac{1}{4}\mE\sup_{s\in[0,t]}|X^{\e,u_{\e}}(s)-X^{0,u_{\e}}(s)|^2+2L_1\e C\mE\int_{0}^{t}\(1+\|X^{0,u_{\e}}_s\|^2_{\infty}+\delta_{X^0_s}(\|\cdot\|^2_{\infty})\)\dif s\\
&&+2L_2C\mE\int_{0}^{t}\(\kappa_1(\|X^{\e,u_{\e}}_s-X^{0,u_\e}_s\|^2_{\infty})+\kappa_2(\mW_2(\sL_{X^{\e}_s},\delta_{X^0_s})^2)\)\dif s.
\de
So, we obtain that
\ce
&&\mE\sup\limits_{s\in[0,t]}|X^{\e,u_{\e}}(s)-X^{0,u_{\e}}(s)|^2\\
&\leq&\frac{1}{2}\mE\sup_{s\in[0,t]}|X^{\e,u_{\e}}(s)-X^{0,u_{\e}}(s)|^2+2L_1\e(C+1)\int_{0}^{t}\(1+\mE\|X^{0,u_{\e}}_s\|^2_{\infty}+\delta_{X^0_s}(\|\cdot\|^2_{\infty})\)\dif s\\
&&+(3+4N+2C)L_2\int_{0}^{t}\(\kappa_1(\mE\|X^{\e,u_{\e}}_s-X^{0,u_\e}_s\|^2_{\infty})+\kappa_2(\mW_2(\sL_{X^{\e}_s},\delta_{X^0_s})^2)\)\dif s.
\de
Note that
\ce
\|X^{0,u_{\e}}_s\|^2_{\infty}=\sup\limits_{\t\in[-r_0,0]}|X^{0,u_{\e}}(s+\t)|^2\leq \sup\limits_{r\in[-r_0,T]}|X^{0,u_{\e}}(r)|^2.
\de
Thus, by Lemma \ref{ldplmlimbnd}, we get that $\|X^{0,u_{\e}}_s\|^2_{\infty}\leq C$, where the constant $C>0$ is independent of $\e$ and $\omega$. Besides, the same deduction to that for (\ref{xex0distes}) yields that
\ce
&&\mE\|X^{\e,u_{\e}}_s-X^{0,u_\e}_s\|^2_{\infty}\leq \mE\sup\limits_{r\in[0,s]}|X^{\e,u_{\e}}(r)-X^{0,u_{\e}}(r)|^2,\\
&&\mW_2(\sL_{X^{\e}_s},\delta_{X^0_s})^2\leq \mE\sup\limits_{r\in[0,s]}|X^{\e}(r)-X^{0}(r)|^2.
\de
Based on the above deduction, we conclude that
\ce
&&\mE\sup\limits_{s\in[0,t]}|X^{\e,u_{\e}}(s)-X^{0,u_{\e}}(s)|^2\\
&\leq& \e C+(3+4N+2C)L_2\kappa_2(\mE\sup\limits_{r\in[0,T]}|X^{\e}(r)-X^{0}(r)|^2)T\\
&&+(3+4N+2C)L_2\int_{0}^{t}\kappa_1(\mE\sup\limits_{r\in[0,s]}|X^{\e,u_{\e}}(r)-X^{0,u_{\e}}(r)|^2)\dif s.
\de
By Lemma \ref{bihari}, it holds that
\ce
\mE\sup\limits_{s\in[0,T]}|X^{\e,u_{\e}}(s)-X^{0,u_{\e}}(s)|^2\leq \Psi^{-1}_1(\Psi_1(C_\e)+(3+4N+2C)L_2T),
\de
where $C_\e:= \e C+(3+4N+2C)L_2\kappa_2(\mE\sup\limits_{r\in[0,T]}|X^{\e}(r)-X^{0}(r)|^2)T$. Taking the limit on two sides of the above inequality, by Lemma \ref{lmcenterdis} and the properties of $\kappa_1$ and $\kappa_2$ we obtain that
\ce
\lim_{\e\to0}\mE\(\sup_{t\in[-r_0,T]}|X^{\e,u_{\e}}(t)-X^{0,u_{\e}}(t)|^2\)=0.
\de

{\bf Step 2.} We prove that for any $\delta>0$, 
$$
\lim_{\e\to0}\mP\left(\sup\limits_{t\in[-r_0,T]}\left|\Gamma^{\e}\left(W+\frac{1}{\sqrt{\e}}\int_{0}^{\cdot}u_{\e}(s)\dif s\right)(t)-\Gamma^{0}\left(\int_{0}^{\cdot}u_{\e}(s)\dif s\right)(t)\right|\geq \delta\right)=0.
$$

By the Chebyshev inequality, it holds that
\ce
\mP\left(\sup\limits_{t\in[-r_0,T]}|X^{\e,u_{\e}}(t)-X^{0,u_{\e}}(t)|\geq \delta\right)\leq \frac{1}{\d^2}\mE\sup_{t\in[-r_0,T]}|X^{\e,u_{\e}}(t)-X^{0,u_{\e}}(t)|^2.
\de
So, by the result in {\bf Step 1.} we get the required convergence. The proof is complete.
\end{proof}

Now, the main result in this section is stated as follows.

\bt\label{ldpth}
Suppose that $({\bf H}_1)$ and $({\bf H}_2)$ hold. Then $\{X^{\e}(\cdot), \e\in(0,1)\}$ satisfies the LDP in $C([-r_0,T], \overline{\cD(A)})$ with the rate function given by 
\ce
I(g)=\frac{1}{2}\inf_{\substack{\{u\in L^2([0,T],\mR^m):\\
g=X^{0,u}\}}}\left\{\int_{0}^{T}|u(s)|^2\dif s\right\}, \quad g\in\mS,
\de
where $\inf\emptyset=+\infty$.
\et
\begin{proof}
To prove Theorem \ref{ldpth}, it is sufficient to verify Condition \ref{cond} in Theorem \ref{ldpcon}. Under the assumptions $({\bf H}_1)$ and $({\bf H}_2)$, $(i)$ and $(ii)$ in Condition \ref{cond} are verified in Proposition \ref{ldpco1} and Proposition \ref{ldpco2}, respectively. The proof is complete.
\end{proof}

\section{The MDP for path-dependent multivalued McKean-Vlasov SDEs}\label{sec:mdp}

In this section, we take $\sC=C([-r_0,0], \overline{\cD(A)})$, require $0\in\cD(A), 0\in A(0)$ and study the MDP for path-dependent multivalued McKean-Vlasov SDEs.

Consider the following path-dependent multivalued McKean-Vlasov SDE:
\be\left\{\begin{array}{ll}
	\dif \frac{\tilde{X}^{\e}(t)-X^0(t)}{a(\e)}\in -A(\frac{\tilde{X}^{\e}(t)-X^0(t)}{a(\e)})\dif t+\frac{b(\tilde{X}^{\e}_t,\sL_{\tilde{X}^{\e}_t})-b(X^0_t,\delta_{X^0_t})}{a(\e)}\dif t+\frac{\sqrt{\e}\s(\tilde{X}^{\e}_t,\sL_{\tilde{X}^{\e}_t})}{a(\e)}\dif W(t), \\
	\qquad\qquad\qquad\qquad t\in(0,T],\\
	\frac{\tilde{X}^{\e}(t)-X^0(t)}{a(\e)}=0,\quad t\in[-r_0,0],
\end{array}
\right.
\label{eqmdp}
\ee
where $X^0_\cdot$ is the solution of Eq.(\ref{eqct}) and $a(\e)$ satisfies 
\be
a(\e)\to0\text{ and } \frac{\e}{a^2(\e)}\to0 \text{ as } \e\to0.
\label{aeps}
\ee

\br
We mention that since $A$ is nonlinear, we can not obtain the equation which $\frac{X^{\e}(t)-X^0(t)}{a(\e)}$ satisfies through the direct computation. So, in order to study the MDP of Eq.(\ref{eq0}), it is necessary to construct Eq.(\ref{eqmdp}).
\er

Under $({\bf H}_1)$ and $({\bf H}_2)$, by Theorem 4.1 in \cite{mq}, Eq.(\ref{eqmdp}) has a unique solution denoted as $(\frac{\tilde{X}^{\e}_\cdot-X^0_\cdot}{a(\e)},\tilde{K}^{\e}(\cdot))$. Let $\tilde{M}^{\e}(\cdot)=\frac{\tilde{X}^{\e}(\cdot)-X^0(\cdot)}{a(\e)}$. To establish the MDP of Eq.(\ref{0eq}), it is equivalent to show that $\tilde{M}^{\e}(\cdot)$ satisfies the LDP. 

Next, in order to prove the LDP for $\tilde{M}^{\e}(\cdot)$ in $\mS=C([-r_0,T], \overline{\cD(A)})$, we give two sufficient conditions.

\bco\label{cond2}
Let $\cG^{\e}: C([0,T], \mR^m)\to \mS$ be a family of measurable mappings. There exists a measurable mapping $\cG^0: C([0,T], \mR^m)\to \mS$ such that
\begin{enumerate}[$(i)$]
\item For each $N\in\mN$, any $\{h_\e\}\subset S^N$ and $h\in S^N$ satisfying $h_\e\to h$ as $\e\rightarrow0$, 
$$
\lim\limits_{\e\to0}\sup\limits_{t\in[-r_0,T]}\left|\cG^0\left(\int_{0}^{\cdot}h_\e(s)\dif s\right)(t)-\cG^{0}\left(\int_{0}^{\cdot}h(s)\dif s\right)(t)\right|=0.
$$
\item For each $N\in\mN$, any $\{u_\e\}\subset \cS^N$ and any $\delta>0$, 	
$$
\lim_{\e\to0}\mP\(\sup_{t\in[-r_0,T]}\left|\cG^{\e}\left(W+\frac{a(\e)}{\sqrt{\e}}\int_{0}^{\cdot}u_{\e}(s)\dif s\right)(t)-\cG^{0}\left(\int_{0}^{\cdot}u_{\e}(s)\dif s\right)(t)\right|\geq \delta\)=0.	
$$
\end{enumerate}
\eco

Under the above assumptions, we have the following result (\cite[Theorem 3.2]{aw}).

\bt\label{mdpcon}
Set $M^{\e}(\cdot):=\cG^{\e}(W)$. Assume that $(i)$ and $(ii)$ in Condition \ref{cond2} hold. Then the family $\{M^{\e}(\cdot), \e\in(0,1)\}$ satisfies the LDP in $\mS$ with the rate function $I$ given by
$$
I(g)=\frac{1}{2}\inf_{\substack{\{u\in L^2([0,T],\mR^m): \\g=\cG^0(\int_{0}^{\cdot}u(s)\dif s)\}}}\left\{\int_{0}^{T}|u(s)|^2\dif s\right\}, \quad g\in \mS,
$$
where $\inf\emptyset=+\infty$.
\et

In order to verify that those mappings associated with $\tilde{M}^{\e}(\cdot)$ satisfy Condition \ref{cond2}, we make the following stronger assumptions:
\begin{enumerate}[$({\bf H}_2^{\prime})$]
\item There exists a constant $L'_2>0$ such that for any $\zeta, \eta\in\sC$, $\mu,\nu\in\sP_2^{\sC}$
\ce
&&|b(\zeta,\mu)-b(\eta,\nu)|^2\leq L_2^{\prime}\left(\|\zeta-\eta\|^2_\infty+\mW_2(\mu,\nu)^2\right),\\
&&\|\s(\zeta,\mu)-\s(\eta,\nu)\|^2\leq L_2^{\prime}\left(\|\zeta-\eta\|^2_\infty+\mW_2(\mu,\nu)^2\right).
\de
\end{enumerate}
\begin{enumerate}[$({\bf H}_3)$] 
\item For any $\mu\in\sP_2^{\sC}$, $b(\cdot,\mu):\sC\to\mR^d$ is Fr\'echet differentiable. There exists a constant $L_3>0$ such that for any $\zeta,\eta\in\sC$, $\mu,\nu\in\sP_2^{\sC}$
\ce
\|Db(\zeta,\mu)-Db(\eta,\nu)\|^2_{\sC^{\ast}}\leq L_3(\|\zeta-\eta\|^2_{\infty}+\mW_2(\mu,\nu)^2),
\de
where $Db(\zeta,\mu)$ stands for the Fr\'echet derivative of $b(\zeta,\mu)$ with respect to $\zeta$, and $\sC^{\ast}$ is the duality space of $\sC$.
\end{enumerate}

\br
$(i)$ It is obvious that $({\bf H}_2^{\prime})$ is stronger than $({\bf H}_2)$.

$(ii)$ $({\bf H}_3)$ implies that for any $\zeta\in\sC$, $\mu\in\sP_2^{\sC}$
\be
\|Db(\zeta,\mu)\|^2_{\sC^{\ast}}&\leq& 2\|Db(0,\d_0)\|^2_{\sC^{\ast}}+2L_3(\|\xi\|^2_\infty+\mu(\|\cdot\|^2_\infty))\no\\
&\leq&C(1+\|\xi\|^2_\infty+\mu(\|\cdot\|^2_\infty)).
\label{dblinegrow}
\ee
\er

Under $({\bf H}_1)$ and $({\bf H}_2^{\prime})$, by the definition of strong solutions for Eq.(\ref{eqmdp}), there exists a unique measurable functional $\cG^{\e}: C([0,T], \mR^m)\to C([-r_0,T], \overline{\cD(A)})$ such that $\tilde{M}^{\e}(\cdot)=\cG^\e(W)$.

Next, consider the following controlled path-dependent multivalued SDE:
\be\left\{\begin{array}{l}
\dif \frac{\tilde{X}^{\e,u}(t)-X^0(t)}{a(\e)} \in -A(\frac{\tilde{X}^{\e,u}(t)-X^0(t)}{a(\e)})\dif t+\frac{b(\tilde{X}^{\e,u}_t,\sL_{\tilde{X}^{\e}_t})-b(X^0_t,\delta_{X^0_t})}{a(\e)}\dif t \\
\qquad\qquad\qquad\quad+\s(\tilde{X}^{\e,u}_t,\sL_{\tilde{X}^{\e}_t})u(t)\dif t+\frac{\sqrt{\e}\s(\tilde{X}^{\e,u}_t,\sL_{\tilde{X}^{\e}_t})}{a(\e)}\dif W(t), \\
\qquad\qquad\qquad\qquad       u\in\cS^N,\quad t\in(0,T],\\
\frac{\tilde{X}^{\e,u}(t)-X^0(t)}{a(\e)}=0,t\in[-r_0,0].
\end{array}
\right.
\label{eqmdpctrl}
\ee
By the Girsanov theorem, Eq.(\ref{eqmdpctrl}) has a unique strong solution denoted as $(\frac{\tilde{X}^{\e,u}_{\cdot}-X^0_{\cdot}}{a(\e)},\tilde{K}^{\e,u}(\cdot))$. Set $\tilde{M}^{\e,u}(\cdot):=\frac{\tilde{X}^{\e,u}(\cdot)-X^0(\cdot)}{a(\e)}$. Moreover, we have $\tilde{M}^{\e,u}(\cdot)=\cG^{\e}(W+\frac{a(\e)}{\sqrt{\e}}\int_{0}^{\cdot}u(s)\dif s)$.
Let $(\tilde{M}^{0,u}_{\cdot},\tilde{K}^{0,u}(\cdot))$ be the solution of the following determinstic equation
\be\left\{\begin{array}{ll}
\dif \tilde{M}^{0,u}(t)\in -A(\tilde{M}^{0,u}(t))\dif t+_{\sC^{\ast}}\<Db(X^0_t,\delta_{X^0_t}),\tilde{M}^{0,u}_t\>_{\sC}\dif t+\s(X^0_t,\delta_{X^0_t})u(t)\dif t, \\
\qquad\qquad\qquad u\in\cS^N, \quad t\in(0,T],\\
\tilde{M}^{0,u}(t)=0,t\in[-r_0,0].
\end{array}
\right.
\label{eqmdplim}
\ee
Define the measurable map $\cG^0: C([0,T], \mR^m)\to C([-r_0,T], \overline{\cD(A)})$ as $\cG^0(\int_{0}^{\cdot}u(s)\dif s)=\tilde{M}^{0,u}(\cdot).$

In the following, we verify $(i)$ and $(ii)$ in Condition \ref{cond2} with $\cG^\e,\cG^0$.

\subsection{Verification of Condition \ref{cond2} $(i)$} 

First, we prepare the following estimate.

\bl\label{mdplimbound}
Assume $({\bf H}_1), ({\bf H}_2^{\prime})$ and $({\bf H}_3)$ hold. Then Eq.(\ref{eqmdplim}) has a unique solution $(\tilde{M}^{0,u}_{\cdot},\tilde{K}^{0,u}(\cdot))$. Moreover, it holds that 
\ce
\sup_{u\in \cS^N}\sup_{t\in[-r_0,T]}|\tilde{M}^{0,u}(t)|^2+\sup_{u\in\cS^N}|\tilde{K}^{0,u}|_0^T\leq C,\quad a.s.
\de
where the constant $C>0$ is independent of $\omega$.
\el
\begin{proof}
Firstly, by the similar proof to that in \cite[Theorem 3.4]{mq}, we obtain that Eq.(\ref{eqmdplim}) has a unique solution $(\tilde{M}^{0,u}_{\cdot},\tilde{K}^{0,u}(\cdot))$. 

Next, for any $v\in \mathrm{Int}(\cD(A))$ and $t\in[0,T]$, it holds that
\ce
|\tilde{M}^{0,u}(t)-v|^2&=&|v|^2-2\int_{0}^{t}\<\tilde{M}^{0,u}(s)-v,\dif \tilde{K}^{0,u}(s)\>\\
&&+2\int_{0}^{t}\<\tilde{M}^{0,u}(s)-v,_{\sC^{\ast}}\<Db(X^0_s,\delta_{X^0_s}),\tilde{M}^{0,u}_s\>_{\sC}\>\dif s\\
&&+2\int_{0}^{t}\<\tilde{M}^{0,u}(s)-v,\s(X^0_s,\delta_{X^0_s})u(s)\>\dif s\\
&=:&|v|^2+I_1(t)+I_2(t)+I_3(t).
\de
By Lemma \ref{L2}, we get that
\ce
I_1(t)&\leq&-2\gamma_1|\tilde{K}^{0,u}|_0^t+2\gamma_2\int_{0}^{t}|\tilde{M}^{0,u}(s)-v|\dif s+2\gamma_3t\\
&\leq&-2\gamma_1|\tilde{K}^{0,u}|_0^t+2\gamma_2\int_{0}^{t}|\tilde{M}^{0,u}(s)-v|^2\dif s+2(\gamma_2+\gamma_3)t.
\de
For $I_2(t)$, by (\ref{dblinegrow}), (\ref{x02moes}) and (\ref{x0d0es}), it holds that
\ce
I_2(t)&\leq&2\int_{0}^{t}|\tilde{M}^{0,u}(s)-v|\|Db(X^0_s,\delta_{X^0_s})\|_{\sC^{\ast}}\|\tilde{M}^{0,u}_s\|_{\infty}\dif s\\
&\leq&\int_{0}^{t}|\tilde{M}^{0,u}(s)-v|^2\dif s+\int_{0}^{t}\|Db(X^0_s,\delta_{X^0_s})\|^2_{\sC^{\ast}}\|\tilde{M}^{0,u}_s\|^2_{\infty}\dif s\\
&\leq&\int_{0}^{t}|\tilde{M}^{0,u}(s)-v|^2\dif s+\int_{0}^{t}C(1+\|X^0_s\|^2_\infty+\delta_{X^0_s}(\|\cdot\|^2_\infty))\|\tilde{M}^{0,u}_s\|^2_{\infty}\dif s\\
&\leq&\int_{0}^{t}|\tilde{M}^{0,u}(s)-v|^2\dif s+C(1+2\sup\limits_{s\in[0,T]}\|X^0_s\|^2_\infty)\int_{0}^{t}\|\tilde{M}^{0,u}_s\|^2_{\infty}\dif s\\
&\leq&\int_{0}^{t}|\tilde{M}^{0,u}(s)-v|^2\dif s+C\int_{0}^{t}\|\tilde{M}^{0,u}_s\|^2_{\infty}\dif s.
\de
Note that
\ce
\|\tilde{M}^{0,u}_s\|^2_{\infty}&=&\sup_{\theta\in[-r_0,0]}|\tilde{M}^{0,u}(s+\theta)|^2\leq\sup_{r\in[-r_0,s]}|\tilde{M}^{0,u}(r)|^2\\
&\leq&2\sup_{r\in[-r_0,s]}|\tilde{M}^{0,u}(r)-v|^2+2|v|^2.
\de
Thus, it follows that
\ce
I_2(t)&\leq& \int_{0}^{t}|\tilde{M}^{0,u}(s)-v|^2\dif s+2C\int_{0}^{t}\sup_{r\in[-r_0,s]}|\tilde{M}^{0,u}(r)-v|^2\dif s+2Ct|v|^2\\
&=&(2C+1)\int_{0}^{t}\sup_{r\in[-r_0,s]}|\tilde{M}^{0,u}(r)-v|^2\dif s+2Ct|v|^2.
\de
As for $I_3(t)$, by $({\bf H}_1)$, (\ref{x02moes}), (\ref{x0d0es}), $u\in\cS^N$, the Young inequality and the H$\ddot{\rm o}$lder inequality, we have
\ce
I_3(t)&\leq&2\int_{0}^{t}|\tilde{M}^{0,u}(s)-v|\|\s(X^0_s,\delta_{X^0_s})\||u(s)|\dif s\\
&\leq& 2\sup_{s\in[0,t]}|\tilde{M}^{0,u}(s)-v|\int_{0}^{t}\|\s(X^0_s,\delta_{X^0_s})\||u(s)|\dif s\\
&\leq&\frac{1}{2}\sup_{s\in[0,t]}|\tilde{M}^{0,u}(s)-v|^2+2\left(\int_{0}^{t}\|\s(X^0_s,\delta_{X^0_s})\|^2\dif s\right)\left(\int_{0}^{t}|u(s)|^2\dif s\right)\\
&\leq&\frac{1}{2}\sup_{s\in[0,t]}|\tilde{M}^{0,u}(s)-v|^2+2NL_1\int_{0}^{t}\(1+\|X^0_s\|^2_{\infty}+\d_{X^0_s}(\|\cdot\|^2_{\infty})\)\dif s\\
&\leq&\frac{1}{2}\sup_{s\in[0,t]}|\tilde{M}^{0,u}(s)-v|^2+2NL_1(1+2\sup\limits_{s\in[0,T]}\|X^0_s\|^2_\infty)t\\
&\leq&\frac{1}{2}\sup_{s\in[0,t]}|\tilde{M}^{0,u}(s)-v|^2+Ct.
\de
So, the above deduction yields that
\ce
|\tilde{M}^{0,u}(t)-v|^2+2\gamma_1|\tilde{K}^{0,u}|_0^t&\leq&|v|^2+(2\gamma_2+2C+1)\int_{0}^{t}\sup_{r\in[-r_0,s]}|\tilde{M}^{0,u}(r)-v|^2\dif s\\
&&+(2\gamma_2+2\gamma_3+2C|v|^2+C)t\\
&&+\frac{1}{2}\sup_{s\in[0,t]}|\tilde{M}^{0,u}(s)-v|^2,
\de
and furthermore,
\be
\sup_{s\in[0,t]}|\tilde{M}^{0,u}(s)-v|^2+4\gamma_1|\tilde{K}^{0,u}|_0^t&\leq&2|v|^2+2(2\gamma_2+2C+1)\int_{0}^{t}\sup_{r\in[-r_0,s]}|\tilde{M}^{0,u}(r)-v|^2\dif s\no\\
&&+2(2\gamma_2+2\gamma_3+2C|v|^2+C)t.
\label{eqmklimbnd}
\ee
Besides, since $\tilde{M}^{0,u}(r)=0$ for $r\in[-r_0,0]$, it holds that
\ce
\sup_{s\in[-r_0,t]}|\tilde{M}^{0,u}(s)-v|^2&\leq&\sup_{s\in[0,t]}|\tilde{M}^{0,u}(s)-v|^2+\sup_{s\in[-r_0,0]}|\tilde{M}^{0,u}(s)-v|^2\\
&=&\sup_{s\in[0,t]}|\tilde{M}^{0,u}(s)-v|^2+|v|^2.
\de
Therefore, combining the above deduction, we conclude that
\ce
\sup_{s\in[-r_0,t]}|\tilde{M}^{0,u}(s)-v|^2&\leq&3|v|^2+2(2\gamma_2+2C+1)\int_{0}^{t}\sup_{r\in[-r_0,s]}|\tilde{M}^{0,u}(r)-v|^2\dif s\\
&&+2(2\gamma_2+2\gamma_3+2C|v|^2+C)t,
\de
which together with the Gronwall inequality yields that
\be
\sup_{s\in[-r_0,T]}|\tilde{M}^{0,u}(s)-v|^2\leq C_1^{\prime \prime}e^{C_2^{\prime \prime}T},
\label{eqmlimbnd}
\ee
where $C_1^{\prime \prime}=3|v|^2+2(2\gamma_2+2\gamma_3+2C|v|^2+C)T$ and $C_2^{\prime \prime}=2(2\gamma_2+2C+1)$.

Finally, by (\ref{eqmlimbnd}), it holds that
\ce
\sup_{s\in[-r_0,T]}|\tilde{M}^{0,u}(s)|^2\leq2C_1^{\prime \prime}e^{C_2^{\prime \prime}T}+2|v|^2.
\de
and further,
\ce
\sup_{u\in \cS^N}\sup_{t\in[-r_0,T]}|\tilde{M}^{0,u}(t)|^2\leq2C_1^{\prime \prime}e^{C_2^{\prime \prime}T}+2|v|^2< \infty.
\de
Moreover, it follows from (\ref{eqmklimbnd}) and (\ref{eqmlimbnd}) that $\sup\limits_{u\in\cS^N}|\tilde{K}^{0,u}|_0^T<\infty$.
\end{proof}

The following result is also needed.

\bl\label{fconv}
Assume that $({\bf H}_1)$ and $({\bf H}^{\prime}_2)$ hold. For each $N\in\mN$, any $\{h_\e\}\subset S^N$ and $h\in S^N$, set
$$
f_\e(t):=\int_{0}^{t}\s(X^{0}_s,\d_{X^0_s})(h_\e(s)-h(s))\dif s, \quad t\in[0,T].
$$
If $h_\e\to h$ in $S^N$ as $\e\rightarrow0$, then $f_\e\to0$ in $C([0,T], \mR^d)$. 
\el

Since the proof of the above lemma is similar to that of Lemma \ref{hehconv}, we omit it.

We now turn to verify Condition \ref{cond2} $(i)$.

\bp\label{mdpco1}
Suppose that $({\bf H}_1),({\bf H}_2^{\prime})$ and $({\bf H}_3)$ hold. Then for each $N\in\mN$, any $\{h_{\e}\}\subset S^N$ and $h\in S^N$ satisfying $h_{\e}\to h$ in $S^N$ as $\e\to0$, it holds that
$$
\lim_{\e\to0}\sup_{t\in[-r_0,T]}\left|\cG^0\(\int_{0}^{\cdot}h_{\e}(s)\dif s\)(t)-\cG^0\(\int_{0}^{\cdot}h(s)\dif s\)(t)\right|=0.
$$
\ep

Since the proof of the above proposition is the same to Proposition \ref{ldpco1}, we also omit it.

\subsection{Verification of Condition \ref{cond2} $(ii)$} Firstly, we prepare three following lemmas.

\bl\label{meconvg}
Under $({\bf H}_1)$ and $({\bf H}_2^{\prime})$ it holds that
\ce
\lim_{\e\to0}\mE\sup_{t\in[-r_0,T]}|\tilde{M}^{\e}(t)|^2=0.
\de
\el
\begin{proof}
For any $t\in[0,T]$, the It\^o formula implies that
\ce
|\tilde{M}^{\e}(t)|^2&=&-2\int_{0}^{t}\<\tilde{M}^{\e}(s),\dif\tilde{K}^{\e}(s)\>+2\int_{0}^{t}\<\tilde{M}^{\e}(s),\frac{b(\tilde{X}^{\e}_s,\sL_{\tilde{X}^{\e}_s})-b(X^0_s,\d_{X^0_s})}{a(\e)}\>\dif s\\
&&+2\int_{0}^{t}\<\tilde{M}^{\e}(s),\frac{\sqrt{\e}\s(\tilde{X}^{\e}_s,\sL_{\tilde{X}^{\e}_s})}{a(\e)}\dif W(s)\>+\frac{\e}{a^2(\e)}\int_{0}^{t}\|\s(\tilde{X}^{\e}_s,\sL_{\tilde{X}^{\e}_s})\|^2\dif s.
\de
By Lemma \ref{L1} and $0\in\cD(A), 0\in A(0)$, it holds that
\ce
\<\tilde{M}^{\e}(s)-0,\dif\tilde{K}^{\e}(s)-0\dif s\>\geq 0,
\de
and
\ce
-2\int_{0}^{t}\<\tilde{M}^{\e}(s),\dif\tilde{K}^{\e}(s)\>=-2\int_{0}^{t}\<\tilde{M}^{\e}(s)-0,\dif\tilde{K}^{\e}(s)-0\dif s\>\leq 0.
\de
So, we infer that
\ce
\mE\sup\limits_{s\in[0,t]}|\tilde{M}^{\e}(s)|^2&\leq& 2\mE\int_{0}^{t}|\tilde{M}^{\e}(s)|\left|\frac{b(\tilde{X}^{\e}_s,\sL_{\tilde{X}^{\e}_s})-b(X^0_s,\d_{X^0_s})}{a(\e)}\right|\dif s\\
&&+2\mE\sup\limits_{s\in[0,t]}\left|\int_{0}^{s}\<\tilde{M}^{\e}(r),\frac{\sqrt{\e}\s(\tilde{X}^{\e}_r,\sL_{\tilde{X}^{\e}_r})}{a(\e)}\dif W(r)\>\right|\\
&&+\frac{\e}{a^2(\e)}\mE\int_{0}^{t}\|\s(\tilde{X}^{\e}_s,\sL_{\tilde{X}^{\e}_s})\|^2\dif s\\
&=:&J_1(t)+J_2(t)+J_3(t).
\de

For $J_1(t)$, $({\bf H}_2^{\prime})$ yields that
\ce
J_1(t)&\leq&\mE\int_{0}^{t}|\tilde{M}^{\e}(s)|^2\dif s+\mE\int_{0}^{t}\frac{\left|b(\tilde{X}^{\e}_s,\sL_{\tilde{X}^{\e}_s})-b(X^0_s,\d_{X^0_s})\right|^2}{a^2(\e)}\dif s\\
&\leq&\mE\int_{0}^{t}\sup\limits_{r\in[-r_0,s]}|\tilde{M}^{\e}(r)|^2\dif s+L_2^\prime\mE\int_{0}^{t}\frac{\|\tilde{X}^{\e}_s-X^0_s\|^2_\infty+\mW_2(\sL_{\tilde{X}^{\e}_s},\d_{X^0_s})^2}{a^2(\e)}\dif s\\
&\leq&\mE\int_{0}^{t}\sup\limits_{r\in[-r_0,s]}|\tilde{M}^{\e}(r)|^2\dif s+2L_2^\prime\int_{0}^{t}\mE\|\tilde{M}^{\e}_s\|^2_\infty\dif s\\
&\leq&(1+2L_2^\prime)\mE\int_{0}^{t}\sup\limits_{r\in[-r_0,s]}|\tilde{M}^{\e}(r)|^2\dif s.
\de
For $J_2(t)$, the Burkholder-Davis-Gundy inequality and $({\bf H}_1)$ imply that
\ce
J_2(t)&\leq& 2C\mE\left(\int_{0}^{t}|\tilde{M}^{\e}(s)|^2\left|\frac{\sqrt{\e}\s(\tilde{X}^{\e}_s,\sL_{\tilde{X}^{\e}_s})}{a(\e)}\right|^2\dif s\right)^{1/2}\\
&\leq& 2C\mE\sup\limits_{s\in[0,t]}|\tilde{M}^{\e}(s)|\left(\int_{0}^{t}\left|\frac{\sqrt{\e}\s(\tilde{X}^{\e}_s,\sL_{\tilde{X}^{\e}_s})}{a(\e)}\right|^2\dif s\right)^{1/2}\\
&\leq&\frac{1}{2}\mE\sup\limits_{s\in[0,t]}|\tilde{M}^{\e}(s)|^2+C\frac{\e}{a^2(\e)}\mE\int_{0}^{t}\|\s(\tilde{X}^{\e}_s,\sL_{\tilde{X}^{\e}_s})\|^2\dif s.
\de
Besides, $({\bf H}_1)$ implies that
\be
\mE\int_{0}^{t}\|\s(\tilde{X}^{\e}_s,\sL_{\tilde{X}^{\e}_s})\|^2\dif s&\leq&L_1\mE\int_{0}^{t}\left(1+\|\tilde{X}^{\e}_s\|^2_{\infty}+\sL_{\tilde{X}^{\e}_s}(\|\cdot\|^2_{\infty})\right)\dif s\no\\
&\leq&L_1\int_{0}^{t}\left(1+2\mE\|a(\e)\tilde{M}^{\e}_s+X^0_s\|^2_{\infty}\right)\dif s\no\\
&\leq&L_1\int_{0}^{t}\left(1+4a^2(\e)\mE\|\tilde{M}^{\e}_s\|^2_{\infty}+4\|X^0_s\|^2_{\infty}\right)\dif s.
\label{s2}
\ee
So, we obtain that 
\ce
J_2(t)&\leq&\frac{1}{2}\mE\sup\limits_{s\in[0,t]}|\tilde{M}^{\e}(s)|^2+4L_1C\e\mE\int_{0}^{t}\sup\limits_{r\in[-r_0,s]}|\tilde{M}^{\e}(r)|^2\dif s\\
&&+CL_1\frac{\e}{a^2(\e)}\mE\int_{0}^{t}(1+4\|X_s^0\|^2_\infty)\dif s.
\de
For $J_3(t)$, by the similar deduction to that for $J_2(t)$, it holds that
\ce
J_3(t)\leq 4L_1\e\mE\int_{0}^{t}\sup\limits_{r\in[-r_0,s]}|\tilde{M}^{\e}(r)|^2\dif s+L_1\frac{\e}{a^2(\e)}\mE\int_{0}^{t}(1+4\|X_s^0\|^2_\infty)\dif s.
\de

Collecting the above deduction, we conclude that
\ce
\mE\sup\limits_{s\in[0,t]}|\tilde{M}^{\e}(s)|^2&\leq& \frac{1}{2}\mE\sup\limits_{s\in[0,t]}|\tilde{M}^{\e}(s)|^2+(1+2L_2^\prime+C\e)\mE\int_{0}^{t}\sup\limits_{r\in[-r_0,s]}|\tilde{M}^{\e}(r)|^2\dif s\\
&&+CL_1\frac{\e}{a^2(\e)}\mE\int_{0}^{t}(1+4\|X_s^0\|^2_\infty)\dif s.
\de
Note that $\tilde{M}^{\e}(s)=0, s\in[-r_0,0]$. Thus, by (\ref{x02moes}) it holds that
\ce
\mE\sup\limits_{s\in[-r_0,t]}|\tilde{M}^{\e}(s)|^2&\leq& 2(1+2L_2^\prime+C\e)\mE\int_{0}^{t}\sup\limits_{r\in[-r_0,s]}|\tilde{M}^{\e}(r)|^2\dif s\\
&&+2CL_1\frac{\e}{a^2(\e)}\mE\int_{0}^{T}(1+4\|X_s^0\|^2_\infty)\dif s\\
&\leq& 2(1+2L_2^\prime+C\e)\mE\int_{0}^{t}\sup\limits_{r\in[-r_0,s]}|\tilde{M}^{\e}(r)|^2\dif s+C\frac{\e}{a^2(\e)}.
\de
The Gronwall inequality yields that
\ce
\mE\sup\limits_{s\in[-r_0,T]}|\tilde{M}^{\e}(s)|^2\leq C\frac{\e}{a^2(\e)}e^{2(1+2L_2^\prime+C\e)T}.
\de
Finally, (\ref{aeps}) implies the required limit. The proof is complete.
\end{proof}

\bl\label{me4bnd}
Assume that $({\bf H}_1)$ and $({\bf H}_2^{\prime})$ hold. Then it holds that
\ce
\mE\(\sup_{t\in[-r_0,T]}|\tilde{M}^{\e}(t)|^4\)\leq C. 
\de
\el
\begin{proof}
For any $v\in \mathrm{Int}(\cD(A))$ and $t\in[0,T]$, by the It\^o formula, it holds that
\ce
&&|\tilde{M}^{\e}(t)-v|^2\\
&=&|v|^2-2\int_{0}^{t}\<\tilde{M}^{\e}(s)-v,\dif\tilde{K}^{\e}(s)\>+2\int_{0}^{t}\<\tilde{M}^{\e}(s)-v,\frac{b(\tilde{X}^{\e}_s,\sL_{\tilde{X}^{\e}_s})-b(X^0_s,\d_{X^0_s})}{a(\e)}\>\dif s\\
&&+2\int_{0}^{t}\<\tilde{M}^{\e}(s)-v,\frac{\sqrt{\e}\s(\tilde{X}^{\e}_s,\sL_{\tilde{X}^{\e}_s})}{a(\e)}\dif W(s)\>+\frac{\e}{a^2(\e)}\int_{0}^{t}\|\s(\tilde{X}^{\e}_s,\sL_{\tilde{X}^{\e}_s})\|^2\dif s.
\de
Lemma \ref{L2} implies that
\ce
-2\int_{0}^{t}\<\tilde{M}^{\e}(s)-v,\dif\tilde{K}^{\e}(s)\>&\leq& -2\gamma_1|\tilde{K}^{\e}|_0^{t}+2\gamma_2\int_{0}^{t}|\tilde{M}^{\e}(s)-v|\dif s+2\gamma_3t\\
&\leq&(2\gamma_2+2\gamma_3)t+2\gamma_2\int_{0}^{t}|\tilde{M}^{\e}(s)-v|^2\dif s.
\de
So, we obtain that
\ce
&&\mE\sup_{s\in[0,t]}|\tilde{M}^{\e}(s)-v|^4\\
&\leq&5(|v|^2+(2\gamma_2+2\gamma_3)t)^2+5(2\gamma_2+1)^2T\int_{0}^{t}\mE|\tilde{M}^{\e}(s)-v|^4\dif s\\
&&+5\mE\sup_{s\in[0,t]}\left(\int_{0}^{s}\left|\frac{b(\tilde{X}^{\e}_r,\sL_{\tilde{X}^{\e}_r})-b(X^0_r,\d_{X^0_r})}{a(\e)}\right|^2\dif r\right)^2\\
&&+20\mE\sup_{s\in[0,t]}\left|\int_{0}^{s}\<\tilde{M}^{\e}(r)-v,\frac{\sqrt{\e}\s(\tilde{X}^{\e}_r,\sL_{\tilde{X}^{\e}_r})}{a(\e)}\dif W(r)\>\right|^2\\
&&+5\mE\sup_{s\in[0,t]}\left|\frac{\e}{a^2(\e)}\int_{0}^{s}\|\s(\tilde{X}^{\e}_r,\sL_{\tilde{X}^{\e}_r})\|^2\dif r\right|^2\\
&=:&5(|v|^2+(2\gamma_2+2\gamma_3)t)^2+5(2\gamma_2+1)^2T\int_{0}^{t}\mE|\tilde{M}^{\e}(s)-v|^4\dif s\\
&&+I_{1}(t)+I_2(t)+I_3(t).
\de

For $I_1(t)$, by $({\bf H}_2^{\prime})$, we deduce that
\ce
I_1(t)&\leq&5T\mE\int_{0}^{t}\left|\frac{b(\tilde{X}^{\e}_s,\sL_{\tilde{X}^{\e}_s})-b(X^0_s,\d_{X^0_s})}{a(\e)}\right|^4\dif s\\
&\leq&5TL_2^{\prime 2}\mE\int_{0}^{t}\frac{\(\|\tilde{X}^{\e}_s-X^0_s\|^2_{\infty}+\mW_2(\sL_{\tilde{X}^{\e}_s},\d_{X^0_s})^2\)^2}{a^4(\e)}\dif s\\
&\leq&5TL_2^{\prime 2}\mE\int_{0}^{t}\frac{\(\|\tilde{X}^{\e}_s-X^0_s\|^2_{\infty}+\mE\|\tilde{X}^{\e}_s-X^0_s\|^2_{\infty}\)^2}{a^4(\e)}\dif s\\
&\leq&20TL_2^{\prime 2}\int_{0}^{t}\mE\|\tilde{M}^{\e}_s\|^4_{\infty}\dif s.
\de
Note that
\be
\|\tilde{M}^{\e}_s\|^4_{\infty}\leq\sup_{r\in[-r_0,s]}|\tilde{M}^{\e}(r)|^4\leq8\sup_{r\in[-r_0,s]}|\tilde{M}^{\e}(r)-v|^4+8|v|^4.
\label{tildmsupa}
\ee
Thus, it holds that
\ce
I_1(t)&\leq&160TL_2^{\prime 2}\int_{0}^{t}\(\mE\sup_{r\in[-r_0,s]}|\tilde{M}^{\e}(r)-v|^4+|v|^4\)\dif s\\
&\leq&160TL_2^{\prime 2}\int_{0}^{t}\mE\sup_{r\in[-r_0,s]}|\tilde{M}^{\e}(r)-v|^4\dif s+160T^2L_2^{\prime 2}|v|^4.
\de

For $I_2(t)$, by the Burkholder-Davis-Gundy inequality, the H\"older inequality and the similar deduction to that for (\ref{s2}), we obtain that
\ce
I_2(t)&\leq&C\mE\int_{0}^{t}|\tilde{M}^{\e}(s)-v|^2\frac{\e}{a^2(\e)}\|\s(\tilde{X}^{\e}_s,\sL_{\tilde{X}^{\e}_s})\|^2\dif s\\
&\leq&C\mE\sup_{s\in[0,t]}|\tilde{M}^{\e}(s)-v|^2\int_{0}^{t}\frac{\e}{a^2(\e)}\|\s(\tilde{X}^{\e}_s,\sL_{\tilde{X}^{\e}_s})\|^2\dif s\\
&\leq&\frac{1}{2}\mE\sup_{s\in[0,t]}|\tilde{M}^{\e}(s)-v|^4+CT\left(\frac{\e}{a^2(\e)}\right)^2\mE\int_{0}^{t}\|\s(\tilde{X}^{\e}_s,\sL_{\tilde{X}^{\e}_s})\|^4\dif s\\
&\leq& \frac{1}{2}\mE\sup_{s\in[0,t]}|\tilde{M}^{\e}(s)-v|^4+\e^2CL^2_1T\mE\int_{0}^{t}\sup_{r\in[-r_0,s]}|\tilde{M}^{\e}(r)-v|^4\dif s\\
&&+CTL^2_1\left(\frac{\e}{a^2(\e)}\right)^2\int_{0}^{t}\left(1+128a^4(\e)|v|^4+16\|X^0_s\|^4_{\infty}\right)\dif s.
\de

In view of $I_3(t)$, by the same deduction to that for $I_2(t)$, it holds that
\ce
I_3(t)&\leq&5T\left(\frac{\e}{a^2(\e)}\right)^2\mE\int_{0}^{t}\|\s(\tilde{X}^{\e}_s,\sL_{\tilde{X}^{\e}_s})\|^4\dif s\\
&\leq&\e^2CL^2_1T\mE\int_{0}^{t}\sup_{r\in[-r_0,s]}|\tilde{M}^{\e}(r)-v|^4\dif s\\
&&+CTL^2_1\left(\frac{\e}{a^2(\e)}\right)^2\int_{0}^{t}\left(1+128a^4(\e)|v|^4+16\|X^0_s\|^4_{\infty}\right)\dif s.
\de

Collecting the above deduction, we obtain that
\ce
&&\mE\sup_{s\in[0,t]}|\tilde{M}^{\e}(s)-v|^4\\
&\leq&5(|v|^2+(2\gamma_2+2\gamma_3)t)^2+160T^2L_2^{\prime 2}|v|^4\\
&&+\(5(2\gamma_2+1)^2T+160TL_2^{\prime 2}+CL^2_1T\e^2\)\mE\int_{0}^{t}\sup_{r\in[-r_0,s]}|\tilde{M}^{\e}(r)-v|^4\dif s\\
&&+CL^2_1T\left(\frac{\e}{a^2(\e)}\right)^2\int_{0}^{t}\left(1+128a^4(\e)|v|^4+16\|X^0_s\|^2_{\infty}\right)\dif s\\
&&+\frac{1}{2}\mE\sup_{s\in[0,t]}|\tilde{M}^{\e}(s)-v|^4.
\de
Since (\ref{aeps}) holds, we may assume that $\max\{\e,a(\e),\frac{\e}{a^2(\e)}\}\leq 1$. Then (\ref{x02moes}) and some calculation imply that
\ce
\mE\sup_{s\in[0,t]}|\tilde{M}^{\e}(s)-v|^4&\leq& C\mE\int_{0}^{t}\sup_{r\in[-r_0,s]}|\tilde{M}^{\e}(r)-v|^4\dif s+C,
\de
and further
\ce
\mE\sup_{s\in[-r_0,t]}|\tilde{M}^{\e}(s)-v|^4&\leq&\mE\sup_{s\in[0,t]}|\tilde{M}^{\e}(s)-v|^4+\mE\sup_{s\in[-r_0,0]}|\tilde{M}^{\e}(s)-v|^4\\
&\leq&C\mE\int_{0}^{t}\sup_{r\in[-r_0,s]}|\tilde{M}^{\e}(r)-v|^4\dif s+C+|v|^4.
\de
The Gronwall inequality yields that
\ce
\mE\sup_{s\in[-r_0,T]}|\tilde{M}^{\e}(s)-v|^4&\leq&(C+|v|^4)e^{CT}.
\de
The proof is complete.
\end{proof}

\bl\label{mdplmctrl4bnd}
Assume that $({\bf H}_1)$ and $({\bf H}_2^{\prime})$ hold. Then for any $u_\e\in \cS^N$ and $u\in \cS^N$, there exists $\e_0>0$ such that for any $0<\e<\e_0$
\ce
\mE\(\sup_{t\in[-r_0,T]}|\tilde{M}^{\e,u_\e}(t)|^4\)\leq C. 
\de
\el
\begin{proof}
For any $v\in \mathrm{Int}(\cD(A))$ and $t\in[0,T]$, by the It\^o formula, it holds that
\ce
&&|\tilde{M}^{\e,u_{\e}}(t)-v|^2\\
&=&|v|^2-2\int_{0}^{t}\<\tilde{M}^{\e,u_{\e}}(s)-v,\dif\tilde{K}^{\e,u_{\e}}(s)\>+2\int_{0}^{t}\<\tilde{M}^{\e,u_{\e}}(s)-v,\frac{b(\tilde{X}^{\e,u_{\e}}_s,\sL_{\tilde{X}^{\e}_s})-b(X^0_s,\d_{X^0_s})}{a(\e)}\>\dif s\\
&&+2\int_{0}^{t}\<\tilde{M}^{\e,u_{\e}}(s)-v,\s(\tilde{X}^{\e,u_{\e}}_s,\sL_{\tilde{X}^{\e}_s})u_{\e}(s)\>\dif s+ 2\int_{0}^{t}\<\tilde{M}^{\e,u_{\e}}(s)-v,\frac{\sqrt{\e}\s(\tilde{X}^{\e,u_{\e}}_s,\sL_{\tilde{X}^{\e}_s})}{a(\e)}\dif W(s)\>\\
&&+\frac{\e}{a^2(\e)}\int_{0}^{t}\|\s(\tilde{X}^{\e,u_{\e}}_s,\sL_{\tilde{X}^{\e}_s})\|^2\dif s.
\de
Besides, Lemma \ref{L2} implies that 
\ce
-2\int_{0}^{t}\<\tilde{M}^{\e,u_{\e}}(s)-v,\dif\tilde{K}^{\e,u_{\e}}(s)\>&\leq& -2\g_1|\tilde{K}^{\e,u_{\e}}|^t_0+2\g_2\int_0^t|\tilde{M}^{\e,u_{\e}}(s)-v|\dif s+2\g_3 t\\
&\leq&2(\gamma_2+\gamma_3)t+2\g_2\int_0^t|\tilde{M}^{\e,u_{\e}}(s)-v|^2\dif s.
\de
Therefore, we get that
\ce
&&\mE\sup_{s\in[0,t]}|\tilde{M}^{\e,u_{\e}}(s)-v|^4\\
&\leq&6(|v|^2+2(\gamma_2+\gamma_3)t)^2+24\gamma_2^2T\mE\int_{0}^{t}|\tilde{M}^{\e,u_{\e}}(s)-v|^4\dif s\\
&&+24\mE\sup_{s\in[0,t]}\left|\int_{0}^{s}\<\tilde{M}^{\e,u_{\e}}(r)-v,\frac{b(\tilde{X}^{\e,u_{\e}}_r,\sL_{\tilde{X}^{\e}_r})-b(X^0_r,\d_{X^0_r})}{a(\e)}\>\dif r\right|^2\\
&&+24\mE\sup_{s\in[0,t]}\left|\int_{0}^{s}\<\tilde{M}^{\e,u_{\e}}(r)-v,\s(\tilde{X}^{\e,u_{\e}}_r,\sL_{\tilde{X}^{\e}_r})u_{\e}(r)\>\dif r\right|^2\\
&&+24\mE\sup_{s\in[0,t]}\left|\int_{0}^{s}\<\tilde{M}^{\e,u_{\e}}(r)-v,\frac{\sqrt{\e}\s(\tilde{X}^{\e,u_{\e}}_r,\sL_{\tilde{X}^{\e}_r})}{a(\e)}\dif W(r)\>\right|^2\\
&&+6\mE\sup_{s\in[0,t]}\left|\frac{\e}{a^2(\e)}\int_{0}^{s}\|\s(\tilde{X}^{\e,u_{\e}}_r,\sL_{\tilde{X}^{\e}_r})\|^2\dif r\right|^2\\
&=:&6\(|v|^2+2(\gamma_2+\gamma_3)t\)^2+24\gamma_2^2T\mE\int_{0}^{t}|\tilde{M}^{\e,u_{\e}}(s)-v|^4\dif s\\
&&+\Sigma_1(t)+\Sigma_2(t)+\Sigma_3(t)+\Sigma_4(t).
\de

For $\Sigma_1(t)$, by $({\bf H}_2^{\prime})$, it holds that
\ce
\Sigma_1(t)&\leq&24T\mE\int_{0}^{t}|\tilde{M}^{\e,u_{\e}}(s)-v|^2\left|\frac{b(\tilde{X}^{\e,u_{\e}}_s,\sL_{\tilde{X}^{\e}_s})-b(X^0_s,\d_{X^0_s})}{a(\e)}\right|^2\dif s\\
&\leq&12T\mE\int_{0}^{t}|\tilde{M}^{\e,u_{\e}}(s)-v|^4\dif s+12T\mE\int_{0}^{t}\left|\frac{b(\tilde{X}^{\e,u_{\e}}_s,\sL_{\tilde{X}^{\e}_s})-b(X^0_s,\d_{X^0_s})}{a(\e)}\right|^4\dif s\\
&\leq&12T\mE\int_{0}^{t}|\tilde{M}^{\e,u_{\e}}(s)-v|^4\dif s+12T{L_2^{\prime}}^2\mE\int_{0}^{t}\(\frac{\|\tilde{X}^{\e,u_{\e}}_s-X^0 _s\|^2_{\infty}+\mW_2(\sL_{\tilde{X}^{\e}_s},\delta_{X^0_s})^2}{a^2(\e)}\)^2\dif s\\
&\leq&12T\mE\int_{0}^{t}|\tilde{M}^{\e,u_{\e}}(s)-v|^4\dif s+24T{L_2^{\prime}}^2\mE\int_{0}^{t}\left(\|\tilde{M}^{\e,u_{\e}}_s\|^4_{\infty}+\mE\|\tilde{M}^{\e}_s\|^4_{\infty}\right)\dif s.
\de
Note that
\ce
\|\tilde{M}^{\e,u_{\e}}_s\|^4_{\infty}\leq\sup_{r\in[-r_0,s]}|\tilde{M}^{\e,u_{\e}}(r)|^4\leq8\sup_{r\in[-r_0,s]}|\tilde{M}^{\e,u_{\e}}(r)-v|^4+8|v|^4.
\de
Thus, we infer that
\ce
\Sigma_1(t)&\leq&12T\mE\int_{0}^{t}|\tilde{M}^{\e,u_{\e}}(s)-v|^4\dif s+24T{L_2^{\prime}}^2\mE\int_{0}^{t}\(8\sup_{r\in[-r_0,s]}|\tilde{M}^{\e,u_{\e}}(r)-v|^4+8|v|^4\)\dif s\\
&&+24T{L_2^{\prime}}^2\mE\int_{0}^{t}\|\tilde{M}^{\e}_s\|^4_{\infty}\dif s\\
&\leq&(12T+192T{L_2^{\prime}}^2)\mE\int_{0}^{t}\sup_{r\in[-r_0,s]}|\tilde{M}^{\e,u_{\e}}(r)-v|^4\dif s+192T^2{L_2^{\prime}}^2|v|^4\\
&&+24T{L_2^{\prime}}^2\mE\int_{0}^{t}\|\tilde{M}^{\e}_s\|^4_{\infty}\dif s.
\de

As for $\Sigma_2(t)$, $({\bf H}_1)$ and $u_{\e}\in \cS^N$ imply that
\ce
&&\Sigma_2(t)\\
&\leq&24T\mE\int_{0}^{t}|\tilde{M}^{\e,u_{\e}}(s)-v|^2\|\s(\tilde{X}^{\e,u_{\e}}_s,\sL_{\tilde{X}^{\e}_s})\|^2|u_{\e}(s)|^2\dif s\\
&\leq&24T\mE\sup_{s\in[0,t]}|\tilde{M}^{\e,u_{\e}}(s)-v|^2\sup_{s\in[0,t]}\|\s(\tilde{X}^{\e,u_{\e}}_s,\sL_{\tilde{X}^{\e}_s})\|^2\int_{0}^{t}|u_{\e}(s)|^2\dif s\\
&\leq&24TN\mE\sup_{s\in[0,t]}|\tilde{M}^{\e,u_{\e}}(s)-v|^2\sup_{s\in[0,t]}\|\s(\tilde{X}^{\e,u_{\e}}_s,\sL_{\tilde{X}^{\e}_s})\|^2\\
&\leq&\frac{1}{8}\mE\sup_{s\in[0,t]}|\tilde{M}^{\e,u_{\e}}(s)-v|^4+C\mE\sup_{s\in[0,t]}\|\s(\tilde{X}^{\e,u_{\e}}_s,\sL_{\tilde{X}^{\e}_s})\|^4\\
&\leq&\frac{1}{8}\mE\sup_{s\in[0,t]}|\tilde{M}^{\e,u_{\e}}(s)-v|^4+CL_1^2\mE\sup_{s\in[0,t]}\left(1+\|\tilde{X}^{\e,u_{\e}}_s\|^2_{\infty}+\sL_{\tilde{X}^{\e}_s}(\|\cdot\|^2_{\infty})\right)^2\\
&\leq&\frac{1}{8}\mE\sup_{s\in[0,t]}|\tilde{M}^{\e,u_{\e}}(s)-v|^4+CL_1^2\mE\sup_{s\in[0,t]}\left(1+\|\tilde{X}^{\e,u_{\e}}_s\|^4_{\infty}+\mE\|\tilde{X}^{\e}_s\|^4_{\infty}\right)\\
&\leq&\frac{1}{8}\mE\sup_{s\in[0,t]}|\tilde{M}^{\e,u_{\e}}(s)-v|^4+CL_1^2\mE\sup_{s\in[0,t]}\left(1+8a^4(\e)\|\tilde{M}^{\e,u_{\e}}_s\|^4_{\infty}+8\|X^0_s\|^4_{\infty}+\mE\|\tilde{X}^{\e}_s\|^4_{\infty}\right)\\
&\leq&\frac{1}{8}\mE\sup_{s\in[0,t]}|\tilde{M}^{\e,u_{\e}}(s)-v|^4+CL_1^2a^4(\e)\mE\sup_{s\in[-r_0,t]}|\tilde{M}^{\e,u_{\e}}(s)-v|^4+CL_1^2a^4(\e)|v|^4\\
&&+CL_1^2\(1+8\sup_{s\in[0,t]}\|X^0_s\|^4_{\infty}+\mE\sup_{s\in[0,t]}\|\tilde{X}^{\e}_s\|^4_{\infty}\).
\de
Since $a(\e)\to0$ as $\e$ goes to 0, we may take $\e_0>0$ such that for any $0<\e<\e_0$, it holds $CL_1^2a^4(\e)\leq\frac{1}{8}$. Thus, we obtain that for any $0<\e<\e_0$,
\ce
\Sigma_2(t)\leq\frac{1}{4}\mE\sup_{s\in[-r_0,t]}|\tilde{M}^{\e,u_{\e}}(s)-v|^4+CL_1^2\(a^4(\e)|v|^4+1+8\sup_{s\in[0,t]}\|X^0_s\|^4_{\infty}+\mE\sup_{s\in[0,t]}\|\tilde{X}^{\e}_s\|^4_{\infty}\).
\de  

For $\Sigma_3(t)$, by the Burkholder-Davis-Gundy inequality and the similar deduction to that for (\ref{s2}), it follows that
\ce
\Sigma_3(t)&\leq& C\mE\int_{0}^{t}|\tilde{M}^{\e,u_{\e}}(s)-v|^2\left\|\frac{\sqrt{\e}\s(\tilde{X}^{\e,u_{\e}}_s,\sL_{\tilde{X}^{\e}_s})}{a(\e)}\right\|^2\dif s\\
&\leq&C\mE\sup_{s\in[0,t]}|\tilde{M}^{\e,u_{\e}}(s)-v|^2\int_{0}^{t}\left\|\frac{\sqrt{\e}\s(\tilde{X}^{\e,u_{\e}}_s,\sL_{\tilde{X}^{\e}_s})}{a(\e)}\right\|^2\dif s\\
&\leq&\frac{1}{4}\mE\sup_{s\in[0,t]}|\tilde{M}^{\e,u_{\e}}(s)-v|^4+CT\left(\frac{\e}{a^2(\e)}\right)^2\mE\int_{0}^{t}\|\s(\tilde{X}^{\e,u_{\e}}_s,\sL_{\tilde{X}^{\e}_s})\|^4\dif s\\
&\leq&\frac{1}{4}\mE\sup_{s\in[0,t]}|\tilde{M}^{\e,u_{\e}}(s)-v|^4+128L_1^2CT\e^2\mE\int_{0}^{t}\sup_{r\in[-r_0,s]}|\tilde{M}^{\e,u_{\e}}(r)-v|^4\dif s\\
&&+4L_1^2CT\left(\frac{\e}{a^2(\e)}\right)^2\int_{0}^{t}\(1+32a^4(\e)|v|^4+4\|X^0_s\|^4_{\infty}+\mE\|\tilde{X}^{\e}_s\|^4_{\infty}\)\dif s.
\de

In view of $\Sigma_4(t)$, by the similar deduction to that for $\Sigma_3(t)$ we have that
\ce
\Sigma_4(t)&\leq& 6T\left(\frac{\e}{a^2(\e)}\right)^2\mE\int_{0}^{t}\|\s(\tilde{X}^{\e,u_{\e}}_s,\sL_{\tilde{X}^{\e}_s})\|^4\dif s\\
&\leq&CL_1^2T\e^2\mE\int_{0}^{t}\sup_{r\in[-r_0,s]}|\tilde{M}^{\e,u_{\e}}(r)-v|^4\dif s\\
&&+24L_1^2T\left(\frac{\e}{a^2(\e)}\right)^2\int_{0}^{t}\(1+32a^4(\e)|v|^4+4\|X^0_s\|^4_{\infty}+\mE\|\tilde{X}^{\e}_s\|^4_{\infty}\)\dif s.
\de
The above deductions yield that
\ce
&&\mE\sup_{s\in[0,t]}|\tilde{M}^{\e,u_{\e}}(s)-v|^4\\
&\leq&6(|v|^2+2(\gamma_2+\gamma_3)t)^2+\frac{1}{2}\mE\sup_{s\in[-r_0,t]}|\tilde{M}^{\e,u_{\e}}(s)-v|^4\\
&&+\left(24\gamma_2^2T+12T+192T{L_2^{\prime}}^2+CL_1^2T\e^2\right)\mE\int_{0}^{t}\sup_{r\in[-r_0,s]}|\tilde{M}^{\e,u_{\e}}(r)-v|^4\dif s\\
&&+192T^2{L_2^{\prime}}^2|v|^4+24T{L_2^{\prime}}^2\mE\int_{0}^{t}\|\tilde{M}^{\e}_s\|^4_{\infty}\dif s\\
&&+CL_1^2\(1+a^4(\e)|v|^4+8\sup_{s\in[0,t]}\|X^0_s\|^4_{\infty}+\mE\sup_{s\in[0,t]}\|\tilde{X}^{\e}_s\|^4_{\infty}\)\\
&&+CL_1^2T\(\frac{\e}{a^2(\e)}\)^2\int_{0}^{t}\(1+32a^4(\e)|v|^4+4\|X^0_s\|^4_{\infty}+\mE\|\tilde{X}^{\e}_s\|^4_{\infty}\)\dif s.
\de
Besides, by (\ref{x02moes}) and Lemma \ref{me4bnd}, we know that
\ce
\sup_{s\in[0,T]}\|X^0_s\|^4_{\infty}\leq C, \quad \mE\sup_{s\in[0,T]}\|\tilde{M}^{\e}_s\|^4_{\infty}\leq C,\quad \mE\sup_{s\in[0,T]}\|\tilde{X}^{\e}_s\|^4_{\infty}\leq C.
\de
So, it holds that
\ce
\mE\sup_{s\in[0,t]}|\tilde{M}^{\e,u_{\e}}(s)-v|^4\leq \frac{1}{2}\mE\sup_{s\in[-r_0,t]}|\tilde{M}^{\e,u_{\e}}(s)-v|^4+C\mE\int_{0}^{t}\sup_{r\in[-r_0,s]}|\tilde{M}^{\e,u_{\e}}(r)-v|^4\dif s+C.
\de
Note that
\ce
\mE\sup_{s\in[-r_0,t]}|\tilde{M}^{\e,u_{\e}}(s)-v|^4&\leq&\mE\sup_{s\in[0,t]}|\tilde{M}^{\e,u_{\e}}(s)-v|^4+\mE\sup_{s\in[-r_0,0]}|\tilde{M}^{\e,u_{\e}}(s)-v|^4\\
&\leq&\mE\sup_{s\in[0,t]}|\tilde{M}^{\e,u_{\e}}(s)-v|^4+|v|^4.
\de
Thus, we conclude that
\ce
\mE\sup_{s\in[-r_0,t]}|\tilde{M}^{\e,u_{\e}}(s)-v|^4\leq 2C\mE\int_{0}^{t}\sup_{r\in[-r_0,s]}|\tilde{M}^{\e,u_{\e}}(r)-v|^4\dif s+2(C+|v|^4),
\de
which together with the Gronwall inequality yields that
\ce
\mE\sup_{s\in[-r_0,T]}|\tilde{M}^{\e,u_{\e}}(s)-v|^4\leq 2(C+|v|^4)e^{2CT}.
\de
The proof is complete.
\end{proof}

At this point, we proceed to verify Condition \ref{cond2} $(ii)$.
\bp\label{mdpco2}
Assume that $({\bf H}_1), ({\bf H}_2^{\prime})$ and $({\bf H}_3)$ hold. Then for each $N\in\mN$, any $\{u_{\e}\}\subset \cS^N$ and any $\d>0$, it holds that
$$
\lim_{\e\to0}\mP\(\sup_{t\in[-r_0,T]}\left|\cG^{\e}\left(W+\frac{a(\e)}{\sqrt{\e}}\int_{0}^{\cdot}u_{\e}(s)\dif s\right)(t)-\cG^{0}\left(\int_{0}^{\cdot}u_{\e}(s)\dif s\right)(t)\right|\geq \d\)=0.	
$$
\ep
\begin{proof}
The proof is divided into two steps. We prove an important limit firstly, and then show the required result.

{\bf Step1.} In this step, we prove that
\ce
\lim_{\e\to0}\mE\(\sup_{t\in[-r_0,T]}\left|\cG^{\e}\left(W+\frac{a(\e)}{\sqrt{\e}}\int_{0}^{\cdot}u_{\e}(s)\dif s\right)(t)-\cG^{0}\left(\int_{0}^{\cdot}u_{\e}(s)\dif s\right)(t)\right|^2\)=0.
\de

By the definition of $\cG^{\e}$ and $\cG^0$, we have
\ce
\tilde{M}^{\e,u_{\e}}(\cdot)&=&\cG^{\e}\left(W+\frac{a(\e)}{\sqrt{\e}}\int_{0}^{\cdot}u_{\e}(s)\dif s\right),\\
\tilde{M}^{0,u_{\e}}(\cdot)&=&\cG^{0}\left(\int_{0}^{\cdot}u_{\e}(s)\dif s\right).
\de
Thus, it suffices to prove
$$
\lim_{\e\to0}\mE\(\sup_{t\in[-r_0,T]}|\tilde{M}^{\e,u_{\e}}(t)-\tilde{M}^{0,u_{\e}}(t)|^2\)=0.
$$

By the It\^o formula, it holds that for any $t\in[0,T]$
\ce
&&|\tilde{M}^{\e,u_{\e}}(t)-\tilde{M}^{0,u_{\e}}(t)|^2\\
&=&-2\int_{0}^{t}\<\tilde{M}^{\e,u_{\e}}(s)-\tilde{M}^{0,u_{\e}}(s),\dif \left(\tilde{K}^{\e,u_{\e}}(s)-\tilde{K}^{0,u_{\e}}(s)\right)\>\\
&&+2\int_{0}^{t}\<\tilde{M}^{\e,u_{\e}}(s)-\tilde{M}^{0,u_{\e}}(s),\frac{b(\tilde{X}^{\e,u_{\e}}_s,\sL_{\tilde{X}^{\e}_s})-b(X^0_s,\d_{X^0_s})}{a(\e)}-_{\sC^{\ast}}\<Db(X^{0}_s,\d_{X^{0}_s}),\tilde{M}^{0,u_\e}_s\>_{\sC}\>\dif s\\
&&+2\int_{0}^{t}\<\tilde{M}^{\e,u_{\e}}(s)-\tilde{M}^{0,u_{\e}}(s),\(\s(\tilde{X}^{\e,u_{\e}}_s,\sL_{\tilde{X}^{\e}_s})-\s(X^{0}_s,\d_{X^0_s})\)u_{\e}(s)\>\dif s\\
&&+2\int_{0}^{t}\<\tilde{M}^{\e,u_{\e}}(s)-\tilde{M}^{0,u_{\e}}(s),\frac{\sqrt{\e}}{a(\e)}\s(\tilde{X}^{\e,u_{\e}}_s,\sL_{\tilde{X}^{\e}_s})\dif W(s)\>\\
&&+\frac{\e}{a^2(\e)}\int_{0}^{t}\|\s(\tilde{X}^{\e,u_{\e}}_s,\sL_{\tilde{X}^{\e}_s})\|^2\dif s\\
&=:&I_1(t)+I_2(t)+I_3(t)+I_4(t)+I_5(t).
\de
According to Lemma \ref{L1}, we know that $I_1(t)\leq0$.

For $I_2(t)$, it holds that
\ce
I_2(t)&=&2\int_{0}^{t}\<\tilde{M}^{\e,u_{\e}}(s)-\tilde{M}^{0,u_{\e}}(s),\frac{b(\tilde{X}^{\e,u_{\e}}_s,\sL_{\tilde{X}^{\e}_s})-b(\tilde{X}^{\e,u_{\e}}_s,\d_{X^0_s})}{a(\e)}\>\dif s\\
&&+2\int_{0}^{t}\<\tilde{M}^{\e,u_{\e}}(s)-\tilde{M}^{0,u_{\e}}(s),\frac{b(\tilde{X}^{\e,u_{\e}}_s,\d_{X^0_s})-b(X^0_s,\d_{X^0_s})}{a(\e)}-_{\sC^{\ast}}\<Db(X^{0}_s,\d_{X^{0}_s}),\tilde{M}^{\e,u_\e}_s\>_{\sC}\>\dif s\\
&&+2\int_{0}^{t}\<\tilde{M}^{\e,u_{\e}}(s)-\tilde{M}^{0,u_{\e}}(s),_{\sC^{\ast}}\<Db(X^{0}_s,\d_{X^{0}_s}),\tilde{M}^{\e,u_\e}_s-\tilde{M}^{0,u_\e}_s\>_{\sC}\>\dif s\\
&=:&I_{21}(t)+I_{22}(t)+I_{23}(t).
\de
Then, by $({\bf H}_2^{\prime})$ and $\tilde{M}^{\e}_{\cdot}=\frac{\tilde{X}^{\e}_\cdot-X^0_\cdot}{a(\e)}$ we have that
\ce
I_{21}(t)&\leq&\int_{0}^{t}|\tilde{M}^{\e,u_{\e}}(s)-\tilde{M}^{0,u_{\e}}(s)|^2\dif s+\int_{0}^{t}\left|\frac{b(\tilde{X}^{\e,u_{\e}}_s,\sL_{\tilde{X}^{\e}_s})-b(\tilde{X}^{\e,u_{\e}}_s,\d_{X^0_s})}{a(\e)}\right|^2\dif s\\
&\leq&\int_{0}^{t}|\tilde{M}^{\e,u_{\e}}(s)-\tilde{M}^{0,u_{\e}}(s)|^2\dif s+L_2^{\prime}\int_{0}^{t}\frac{\mW_2(\sL_{\tilde{X}^{\e}_s},\d_{X^0_s})^2}{a^2(\e)}\dif s\\
&\leq&\int_{0}^{t}|\tilde{M}^{\e,u_{\e}}(s)-\tilde{M}^{0,u_{\e}}(s)|^2\dif s+L_2^{\prime}\int_{0}^{t}\mE\|\tilde{M}^{\e}_s\|^2_{\infty}\dif s.
\de
And $({\bf H}_3)$ and $\tilde{M}^{\e,u}_\cdot=\frac{\tilde{X}^{\e,u}_\cdot-X^0_\cdot}{a(\e)}$ imply that
\ce
I_{22}(t)&\leq&\int_{0}^{t}\left|\frac{b(a(\e)\tilde{M}^{\e,u_{\e}}_s+X^0_s,\d_{X^0_s})-b(X^0_s,\d_{X^{0}_s})}{a(\e)}-_{\sC^{\ast}}\<Db(X^{0}_s,\d_{X^{0}_s}),\tilde{M}^{\e,u_\e}_s\>_{\sC}\right|^2\dif s\\
&&+\int_{0}^{t}|\tilde{M}^{\e,u_{\e}}(s)-\tilde{M}^{0,u_{\e}}(s)|^2\dif s\\
&\leq&\int_{0}^{t}\left|\int_{0}^{1}{_{\sC^{\ast}}}\<Db(r a(\e)\tilde{M}^{\e,u_{\e}}_s+X^0_s,\d_{X^0_s}),\tilde{M}^{\e,u_{\e}}_s\>_{\sC}\dif r-_{\sC^{\ast}}\<Db(X^{0}_s,\d_{X^{0}_s}),\tilde{M}^{\e,u_\e}_s\>_{\sC}\right|^2\dif s\\
&&+\int_{0}^{t}|\tilde{M}^{\e,u_{\e}}(s)-\tilde{M}^{0,u_{\e}}(s)|^2\dif s\\
&\leq&\int_{0}^{t}\[\int_{0}^{1}\|Db(r a(\e)\tilde{M}^{\e,u_{\e}}_s+X^0_s,\d_{X^0_s})-Db(X^{0}_s,\d_{X^{0}_s})\|_{\sC^{\ast}}\dif r\]^2\|\tilde{M}^{\e,u_{\e}}_s\|^2_{\infty}\dif s\\
&&+\int_{0}^{t}|\tilde{M}^{\e,u_{\e}}(s)-\tilde{M}^{0,u_{\e}}(s)|^2\dif s\\
&\leq&\int_{0}^{t}\[\int_{0}^{1}L^{1/2}_3 r a(\e)\|\tilde{M}^{\e,u_{\e}}_s\|_{\infty}\dif r\]^2 \|\tilde{M}^{\e,u_{\e}}_s\|^2_{\infty}\dif s+\int_{0}^{t}|\tilde{M}^{\e,u_{\e}}(s)-\tilde{M}^{0,u_{\e}}(s)|^2\dif s\\
&\leq&L_3 a^2(\e)\int_{0}^{t} \|\tilde{M}^{\e,u_{\e}}_s\|^4_{\infty}\dif s+\int_{0}^{t}|\tilde{M}^{\e,u_{\e}}(s)-\tilde{M}^{0,u_{\e}}(s)|^2\dif s.
\de
Moreover, from (\ref{dblinegrow}), (\ref{x02moes}) and (\ref{x0d0es}) it follows that
\ce
I_{23}(t)&\leq&\int_{0}^{t}\|Db(X^{0}_s,\d_{X^{0}_s})\|^2_{\sC^{\ast}}\|\tilde{M}^{\e,u_{\e}}_s-\tilde{M}^{0,u_{\e}}_s\|^2_{\infty}\dif s+\int_{0}^{t}|\tilde{M}^{\e,u_{\e}}(s)-\tilde{M}^{0,u_{\e}}(s)|^2\dif s\\
&\leq&C(1+2\sup\limits_{s\in[0,T]}\|X^0_s\|^2_\infty)\int_{0}^{t}\sup_{r\in[-r_0,s]}|\tilde{M}^{\e,u_{\e}}(r)-\tilde{M}^{0,u_{\e}}(r)|^2\dif s\\
&&+\int_{0}^{t}|\tilde{M}^{\e,u_{\e}}(s)-\tilde{M}^{0,u_{\e}}(s)|^2\dif s\\
&\leq&(C+1)\int_{0}^{t}\sup_{r\in[-r_0,s]}|\tilde{M}^{\e,u_{\e}}(r)-\tilde{M}^{0,u_{\e}}(r)|^2\dif s.
\de
Therefore, it holds that
\ce
I_2(t)&\leq&(C+3)\int_{0}^{t}\sup_{r\in[-r_0,s]}|\tilde{M}^{\e,u_{\e}}(r)-\tilde{M}^{0,u_{\e}}(r)|^2\dif s\\
&&+L_3a^2(\e)\int_{0}^{t} \|\tilde{M}^{\e,u_{\e}}_s\|^4_{\infty}\dif s+L_2^{\prime}\int_{0}^{t}\mE\|\tilde{M}^{\e}_s\|^2_{\infty}\dif s.
\de

In view of $I_3(t)$, by the H$\ddot{\rm o}$lder inequality, $({\bf H}_2^{\prime})$ and $u_{\e}\in \cS^N$, we get
\ce
I_3(t)&\leq&2\int_{0}^{t}|\tilde{M}^{\e,u_{\e}}(s)-\tilde{M}^{0,u_{\e}}(s)|\left |\left(\s(\tilde{X}^{\e,u_{\e}}_s,\sL_{\tilde{X}^{\e}_s})-\s(X^{0}_s,\d_{X^0_s})\right)u_{\e}(s)\right|\dif s\\
&\leq&2\sup_{s\in[0,t]}|\tilde{M}^{\e,u_{\e}}(s)-\tilde{M}^{0,u_{\e}}(s)|\int_{0}^{t}\left |\left(\s(\tilde{X}^{\e,u_{\e}}_s,\sL_{\tilde{X}^{\e}_s})-\s(X^{0}_s,\d_{X^0_s})\right)u_{\e}(s)\right|\dif s\\
&\leq&\frac{1}{4}\sup_{s\in[0,t]}|\tilde{M}^{\e,u_{\e}}(s)-\tilde{M}^{0,u_{\e}}(s)|^2\\
&&+C\[\int_{0}^{t}\|\s(a(\e)\tilde{M}^{\e,u_{\e}}_s+X^0_s,\sL_{\tilde{X}^{\e}_s})-\s(X^{0}_s,\d_{X^0_s})\|^2\dif s\]\[\int_{0}^{t}|u_{\e}(s)|^2\dif s\]\\ 
&\leq&\frac{1}{4}\sup_{s\in[0,t]}|\tilde{M}^{\e,u_{\e}}(s)-\tilde{M}^{0,u_{\e}}(s)|^2+CNL_2^{\prime}\int_{0}^{t}\(\|a(\e)\tilde{M}^{\e,u_{\e}}_s\|^2_{\infty}+\mW_2(\sL_{\tilde{X}^{\e}_s},\delta_{X^0_s})^2\)\dif s\\
&=&\frac{1}{4}\sup_{s\in[0,t]}|\tilde{M}^{\e,u_{\e}}(s)-\tilde{M}^{0,u_{\e}}(s)|^2+a^2(\e)CNL_2^{\prime}\int_{0}^{t}\(\|\tilde{M}^{\e,u_{\e}}_s\|^2_{\infty}+\mE\|\tilde{M}^{\e}_s\|^2_{\infty}\)\dif s.
\de

For $I_5$, by the similar deduction to that for (\ref{s2}), it holds that
\ce
I_5(t)\leq2\e L_1\int_{0}^{t}\left(\|\tilde{M}^{\e,u_{\e}}_s\|^2_{\infty}+\mE\|\tilde{M}^{\e}_s\|^2_{\infty}\right)\dif s+\frac{\e}{a^2(\e)}L_1\int_{0}^{t}\left(1+4\|X^0_s\|^2_{\infty}\right)\dif s.
\de

Combining the above estimates, we derive that
\ce
&&\sup_{s\in[0,t]}|\tilde{M}^{\e,u_{\e}}(s)-\tilde{M}^{0,u_{\e}}(s)|^2\\
&\leq&\frac{1}{4}\sup_{s\in[0,t]}|\tilde{M}^{\e,u_{\e}}(s)-\tilde{M}^{0,u_{\e}}(s)|^2+(C+3)\int_{0}^{t}\sup_{r\in[-r_0,s]}|\tilde{M}^{\e,u_{\e}}(r)-\tilde{M}^{0,u_{\e}}(r)|^2\dif s\\
&&+L_3a^2(\e)\int_{0}^{t} \|\tilde{M}^{\e,u_{\e}}_s\|^4_{\infty}\dif s+L_2^{\prime}\int_{0}^{t}\mE\|\tilde{M}^{\e}_s\|^2_{\infty}\dif s\\
&&+\left(a^2(\e)CNL_2^{\prime}+2\e L_1\right)\int_{0}^{t}\left(\|\tilde{M}^{\e,u_{\e}}_s\|^2_{\infty}+\mE\|\tilde{M}^{\e}_s\|^2_{\infty}\right)\dif s\\
&&+\frac{\e}{a^2(\e)}L_1\int_{0}^{t}\left(1+4\|X^0_s\|^2_{\infty}\right)\dif s+\sup_{s\in[0,t]}|I_4(s)|.
\de

Next, from the Burkholder-Davis-Gundy inequality and the same deduction to that of $I_5$, it follows that
\ce
&&\mE\sup_{s\in[0,t]}|I_4(s)|\\
&\leq&C\frac{\sqrt{\e}}{a(\e)}\mE\(\int_{0}^{t}|\tilde{M}^{\e,u_{\e}}(s)-\tilde{M}^{0,u_{\e}}(s)|^2\|\s(\tilde{X}^{\e,u_{\e}}_s,\sL_{\tilde{X}^{\e}_s})\|^2\dif s\)^{1/2}\\
&\leq&C\frac{\sqrt{\e}}{a(\e)}\mE\sup_{s\in[0,t]}|\tilde{M}^{\e,u_{\e}}(s)-\tilde{M}^{0,u_{\e}}(s)|\left(\int_{0}^{t}\|\s(\tilde{X}^{\e,u_{\e}}_s,\sL_{\tilde{X}^{\e}_s})\|^2\dif s\right)^{1/2}\\
&\leq&\frac{1}{4}\mE\sup_{s\in[0,t]}|\tilde{M}^{\e,u_{\e}}(s)-\tilde{M}^{0,u_{\e}}(s)|^2+C\frac{\e}{a^2(\e)}\mE\int_{0}^{t}\|\s(\tilde{X}^{\e,u_{\e}}_s,\sL_{\tilde{X}^{\e}_s})\|^2\dif s\\
&\leq&\frac{1}{4}\mE\sup_{s\in[0,t]}|\tilde{M}^{\e,u_{\e}}(s)-\tilde{M}^{0,u_{\e}}(s)|^2+2\e CL_1\mE\int_{0}^{t}\left(\|\tilde{M}^{\e,u_{\e}}_s\|^2_{\infty}+\|\tilde{M}^{\e}_s\|^2_{\infty}\right)\dif s\\
&&+\frac{\e}{a^2(\e)}CL_1\int_{0}^{t}\left(1+4\|X^0_s\|^2_{\infty}\right)\dif s.
\de
So, we obtain that
\ce
&&\mE\sup_{s\in[0,t]}|\tilde{M}^{\e,u_{\e}}(s)-\tilde{M}^{0,u_{\e}}(s)|^2\\
&\leq&\frac{1}{2}\mE\sup_{s\in[0,t]}|\tilde{M}^{\e,u_{\e}}(s)-\tilde{M}^{0,u_{\e}}(s)|^2+(C+3)\mE\int_{0}^{t}\sup_{r\in[-r_0,s]}|\tilde{M}^{\e,u_{\e}}(r)-\tilde{M}^{0,u_{\e}}(r)|^2\dif s\\
&&+L_3a^2(\e)\mE\int_{0}^{t} \|\tilde{M}^{\e,u_{\e}}_s\|^4_{\infty}\dif s+\(L_2^{\prime}+a^2(\e)CNL_2^{\prime}+2\e CL_1\)\int_{0}^{t}\mE\|\tilde{M}^{\e}_s\|^2_{\infty}\dif s\\
&&+\left(a^2(\e)CNL_2^{\prime}+2\e CL_1\right)\mE\int_{0}^{t}\|\tilde{M}^{\e,u_{\e}}_s\|^2_{\infty}\dif s+\frac{\e}{a^2(\e)}CL_1\int_{0}^{t}\left(1+4\|X^0_s\|^2_{\infty}\right)\dif s.
\de
By (\ref{x02moes}), we know $\sup\limits_{s\in[0,T]}\|X^0_s\|^2_{\infty}\leq C$. And by Lemma \ref{mdplmctrl4bnd}, there exists a $\e_0>0$ such that for any $\e\in(0,\e_0)$, it holds that 
$$
\mE\sup\limits_{s\in[0,T]}\|\tilde{M}^{\e,u_{\e}}_s\|^4_{\infty}\leq C, \quad \mE\sup\limits_{s\in[0,T]}\|\tilde{M}^{\e,u_{\e}}_s\|^2_{\infty}\leq C. 
$$
Note that $\tilde{M}^{\e,u_\e}(s)=\tilde{M}^{0,u_\e}(s)=0,s\in[-r_0,0]$. Thus, for any $\e\in(0,\e_0)$, we conclude that
\ce
&&\mE\sup_{s\in[-r_0,t]}|\tilde{M}^{\e,u_{\e}}(s)-\tilde{M}^{0,u_{\e}}(s)|^2\\
&\leq&2(C+3)\mE\int_{0}^{t}\sup_{r\in[-r_0,s]}|\tilde{M}^{\e,u_{\e}}(r)-\tilde{M}^{0,u_{\e}}(r)|^2\dif s+C\left(\e+a^2(\e)+\frac{\e}{a^2(\e)}\right)\\ 
&&+2\(L_2^{\prime}+a^2(\e)CNL_2^{\prime}+2\e CL_1\)T\mE\sup\limits_{s\in[0,T]}\|\tilde{M}^{\e}_s\|^2_{\infty},
\de
which together with the Gronwall inequality yields that
\ce
&&\mE\sup_{s\in[-r_0,T]}|\tilde{M}^{\e,u_{\e}}(s)-\tilde{M}^{0,u_{\e}}(s)|^2\\
&\leq&\left[C\left(\e+a^2(\e)+\frac{\e}{a^2(\e)}\right)+2\(L_2^{\prime}+a^2(\e)CNL_2^{\prime}+2\e CL_1\)T\mE\sup\limits_{s\in[0,T]}\|\tilde{M}^{\e}_s\|^2_{\infty}\right]e^{2(C+3)T}.
\de
Taking the limit on both sides of the above inequality, by (\ref{aeps}) and Lemma \ref{meconvg} we derive that
\ce
\lim_{\e\to0}\mE\(\sup_{t\in[-r_0,T]}|\tilde{M}^{\e,u_{\e}}(t)-\tilde{M}^{0,u_{\e}}(t)|^2\)=0.
\de

{\bf Step 2.} We prove that for any $\d>0$, it holds that
\ce
\lim_{\e\to0}\mP\left(\sup_{t\in[-r_0,T]}\left|\cG^{\e}\left(W+\frac{a(\e)}{\sqrt{\e}}\int_{0}^{\cdot}u_{\e}(s)\dif s\right)(t)-\cG^{0}\left(\int_{0}^{\cdot}u_{\e}(s)\dif s\right)(t)\right|\geq \d\right)=0.
\de

By the Chebyshev inequality, we have that
\ce
\mP\left(\sup_{t\in[-r_0,T]}|\tilde{M}^{\e,u_{\e}}(t)-\tilde{M}^{0,u_{\e}}(t)|^2\geq\d\right)\leq\frac{1}{\d^2}\mE\sup_{t\in[-r_0,T]}|\tilde{M}^{\e,u_{\e}}(t)-\tilde{M}^{0,u_{\e}}(t)|^2,
\de
which together with the result in {\bf Step 1.} yields the desired convergence. The proof is complete.
\end{proof}

The main result of this section is stated as follows.

\bt\label{mdpth01}
Suppose that $({\bf H}_1), ({\bf H}_2^{\prime})$ and $({\bf H}_3)$ hold. Then $\{\tilde{M}^{\e}, \e\in(0,1)\}$ satisfies the LDP on $C([-r_0,T],\overline{\cD(A)})$ with the rate function given by
$$
I(g)=\frac{1}{2}\inf_{\substack{\{u\in L^2([0,T];\mR^m): \\ g=\tilde{M}^{0,u}\}}}\left\{\int_{0}^{T}|u(s)|^2\dif s\right\},
$$
where $\inf\emptyset=+\infty$.
\et
\begin{proof}
By Theorem \ref{mdpcon}, it suffices to verify $(i)$ and $(ii)$ in Condition \ref{cond2}. Under the assumptions $({\bf H}_1), ({\bf H}_2^{\prime})$ and $({\bf H}_3)$, $(i)$ and $(ii)$ are justified in Proposition \ref{mdpco1} and Proposition \ref{mdpco2}, respectively. The proof is complete.
\end{proof}

\section{The CLT for path-dependent multivalued McKean-Vlasov SDEs}\label{sec:clt}

In this section, we take $\sC=C([-r_0,0], \overline{\cD(A)})$, require $0\in\overline{\cD(A)}$, and study the CLT for path-dependent multivalued McKean-Vlasov SDEs.

\subsection{$L$-derivatives for the functions on $\sP_{2}^{\sC}$}

In this subsection, we introduce the $L$-derivatives for the functions on $\sP^\sC_2$ (\cite{brw}). 

For any $\mu\in \sP^\sC_2$, the tangent space at $\mu$ is given by
\ce
T_{\mu,2}:=\{\phi:\sC\to\sC, \phi \text{ is measurable with } \mu(\|\phi(\cdot)\|^2_{\infty})<\infty\},
\de
with the norm defined by $\|\phi\|_{T_{\mu,2}}:=\(\mu(\|\phi(\cdot)\|^2_{\infty})\)^{1/2}$, for $\phi\in T_{\mu,2}$.

\bd
Let $f:\sP^\sC_2\to \mR$ be a continuous function and $\mathrm{Id}$ be the identity map on $\sC$.
\begin{enumerate}[(i)]
\item If for any $\mu\in\sP^{\sC}_2$,
\ce
T_{\mu,2}\ni\phi\to D^{L}_{\phi}f(\mu):=\lim_{\e\to0}\frac{f(\mu\circ(\mathrm{Id}+\e\phi)^{-1})-f(\mu)}{\e}\in \mR
\de
is a well-defined bounded linear functional, then $f$ is called  intrinsically differentiable at $\mu$ and $D^L_{\phi}f(\mu)$ is called the intrinsic derivative of $f$ at $\mu$. In this case, there exists the unique element $D^Lf(\mu)\in T_{\mu,2}^*$ such that 
\ce
{_{T_{\mu,2}^*}}\<D^Lf(\mu),\phi\>_{T_{\mu,2}}:=\int_{\sC}{_{\sC^{\ast}}}\<D^Lf(\mu)(x),\phi(x)\>_{\sC}\mu(\dif x)=D^L_{\phi}f(\mu), \quad \phi\in T_{\mu,2}.
\de
\item If 
\ce
\lim_{\|\phi\|_{T_{\mu,2}}\to0}\frac{|f(\mu\circ(\mathrm{Id}+\phi)^{-1})-f(\mu)-D^L_{\phi}f(\mu)|}{\|\phi\|_{T_{\mu,2}}}=0,
\de
then $f$ is called $L$-differentiable at $\mu$ with the $L$-derivative (i.e. Lions derivative) $D^Lf(\mu)$.
\end{enumerate}
\ed 

Let $C^1(\sP^{\sC}_2)$ denote the set of all continuous functions $f$ on $\sP^{\sC}_2$ satisfying that $f$ is $L$-differentiable at any point $\mu\in \sP^{\sC}_2$, and the $L$-derivative $D^Lf(\mu)(x)$ has a version jointly continuous in $(\mu,x)\in \sP^{\sC}_2\times \sC$. Moreover, let $C^1_b(\sP^{\sC}_2)$ denote the subset of $C^1(\sP^{\sC}_2)$ with $D^Lf(\mu)(x)$ bounded.

For a vector-valued function $f=(f_i)$ or a matrix-valued function $f=(f_{ij})$ with $L$-differentiable components, we write
\ce
D^L_{\phi}f(\mu)=(D^L_{\phi}f_i(\mu)) \text{ or }  D^L_{\phi}f(\mu)=(D^L_{\phi}f_{ij}(\mu)),\ \mu\in\sP^{\sC}_2.
\de 

The following result is from \cite[Theorem 2.1]{brw}.

\bl
Let $f:\sP^{\sC}_2\to \mR$ be a continuous function and $(\varsigma_{\e})_{\e\in[0,1]}$ be a family of $\sC$-valued random variables on a complete probability space $(\Omega,\sF,\mP)$ such that $\dot{\varsigma}_0:=\lim\limits_{\e\to0}\frac{\varsigma_{\e}-\varsigma_0}{\e}$ exists in $L^2(\Omega)$. Let $\mu_0=\sL_{\varsigma_0}$ be atomless. If $f$ is $L$-differentiable such that $D^Lf(\mu_0)$ has a continuous version satisfying
\ce
\|D^Lf(\mu_0)(x)\|_{\sC^{\ast}}\leq C(1+\|x\|_{\sC}),\ x\in\sC
\de
for some constant $C>0$, then
\ce
\lim_{\e\to0}\frac{f(\sL_{\varsigma_\e})-f(\sL_{\varsigma_0})}{\e}=\mE[_{\sC^{\ast}}\<D^Lf(\mu_0)(\varsigma_0),\dot{\varsigma_0}\>_{\sC}].
\de
\el

\subsection{The CLT}

In this subsection, we first prove the CLT of path-dependent multivalued McKean-Vlasov SDEs and then give an example to explain our result.

First of all, since the operator $A$ is nonlinear, in order to study the CLT of Eq.(\ref{eq0}) we construct the following path-dependent multivalued McKean-Vlasov SDE:
\be\left\{\begin{array}{ll}
\dif \frac{\hat{X}^{\e}(t)-X^0(t)}{\sqrt{\e}}\in -A(\frac{\hat{X}^{\e}(t)-X^0(t)}{\sqrt{\e}})\dif t+\frac{b(\hat{X}^{\e}_t,\sL_{\hat{X}^{\e}_t})-b(X^0_t,\delta_{X^0_t})}{\sqrt{\e}}\dif t+\s(\hat{X}^{\e}_t,\sL_{\hat{X}^{\e}_t})\dif W(t), \\
\qquad\qquad\qquad\qquad t\in(0,T],\\
\frac{\hat{X}^{\e}(t)-X^0(t)}{\sqrt{\e}}=0,\quad t\in[-r_0,0].
\end{array}
\right.
\label{eqcltdistb}
\ee
Under $({\bf H}_1)$ and $({\bf H}_2)$, by Theorem 4.1 in \cite{mq}, Eq.(\ref{eqcltdistb}) has a unique solution denoted as $(\frac{\hat{X}^{\e}_\cdot-X^0_\cdot}{\sqrt{\e}},\hat{K}^{\e}(\cdot))$. Set $Z^{\e}(\cdot):=\frac{\hat{X}^{\e}(\cdot)-X^0(\cdot)}{\sqrt{\e}}$. 

Assume:
\begin{enumerate}[$({\bf H}_4)$]  
\item  For $\zeta\in\sC$, $b(\zeta, \cdot)$ is $L$-differentiable. Moreover, there exists a constant $L_4>0$ such that for any $\zeta\in\sC$ and $X,Y,\phi\in L^2(\Omega\to\sC,\mP)$
\ce
&&\big|\mE[_{\sC^{\ast}}\<D^{L}b(\zeta,\sL_X)(X),\phi\>_{\sC}]-\mE[_{\sC^{\ast}}\<D^{L}b(\zeta,\sL_Y)(Y),\phi\>_{\sC}]\big|\\
&&\leq L_4\(\mW_2(\sL_X,\sL_Y)+(\mE\|X-Y\|^2_{\infty})^{1/2}\)\(\mE\|\phi\|^2_{\infty}\)^{1/2}.
\de
And for any $(\zeta,\mu)\in\sC\times\sP^{\sC}_2$,
\ce
\|D^Lb(\zeta,\mu)\|_{T^*_{\mu,2}}\leq L_4.
\de
\end{enumerate}

Let us consider the following path-dependent multivalued McKean-Vlasov SDE:
\be\left\{\begin{array}{ll}
\dif Z(t)\in -A(Z(t))\dif t+_{\sC^{\ast}}\<Db(X^0_t,\d_{X^0_t}),Z_t\>_{\sC}\dif t+\mE[_{\sC^{\ast}}\<D^Lb(X^0_t,\d_{X^0_t})(X^0_t),Z_t\>_{\sC}]\dif t\\
\quad\quad\quad\quad+\s(X^0_t,\d_{X^0_t})\dif W(t), \quad t\in(0,T],\\
Z(t)=0,\quad t\in[-r_0,0].
\end{array}
\right.
\label{eqcltlimit}
\ee
Note that the coefficients of Eq.(\ref{eqcltlimit}) are Lipschitz continuous. Thus, Theorem 4.1 in \cite{mq} implies that Eq.(\ref{eqcltlimit}) has a unique solution $(Z_{\cdot},\hat{K}^0(\cdot))$. Our aim is to show that
\ce
Z^{\e}(\cdot)\xrightarrow{d} Z(\cdot).
\de

We prepare some estimates.
\bl\label{lemma:clt-X0bnd}
Assume $({\bf H}_1)$ and $({\bf H}_2)$ hold. It holds that for any $p\geq1$,
\ce
\sup_{t\in[-r_0,T]}|X^0(t)|^{4p}\leq C.
\de
\el

\bl\label{lemma:clt-Zbnd}
Assume that $({\bf H}_1), ({\bf H}^{\prime}_2), ({\bf H}_3)$ and $({\bf H}_4)$ hold. Then it holds that for any $p\geq1$,
\ce
\mE\left(\sup_{t\in[-r_0,T]}|Z^\e(t)|^{4p}\right)\leq C,\quad \mE\left(\sup_{t\in[-r_0,T]}|Z(t)|^{4p}\right)\leq C.
\de
\el

Since the proofs of two above lemmas are similar to that for (\ref{x02moes}) and Lemma \ref{me4bnd} respectively, we omit them.

Finally, we present the main result of this section and prove it.

\bt\label{thm:clt}
Under assumptions $({\bf H}_1), ({\bf H}^{\prime}_2), ({\bf H}_3)$ and $({\bf H}_4)$, it holds that for any $p\geq1$,
\ce
\mE\left(\sup_{t\in[-r_0,T]}\left|Z^\e(t)-Z(t)\right|^{2p}\right)\leq C\e^p,
\de
where the constant $C>0$ is independent of $\e$.
\et
\begin{proof} 
By the It\^o formula, it holds that
\ce
&&|Z^\e(t)-Z(t)|^{2p}\\
&=&-2p\int_0^t|Z^\e(s)-Z(s)|^{2p-2}\Bigl<Z^\e(s)-Z(s),\dif \hat{K}^\e(s)-\dif \hat{K}^0(s)\Bigr>\\
&&+2p\int_0^t|Z^\e(s)-Z(s)|^{2p-2}\Bigl<Z^\e(s)-Z(s),\frac{b(\hat{X}^\e_s,\sL_{\hat{X}^\e_s})-b(X^0_s,\d_{X^0_s})}{\sqrt{\e}}\\
&&\qquad\qquad -{_{\sC^{\ast}}}\<Db(X^0_s,\d_{X^0_s}),Z_s\>_{\sC}-\mE[_{\sC^{\ast}}\<D^Lb(X^0_s,\d_{X^0_s})(X^0_s),Z_s\>_{\sC}]\Bigr>\dif s\\
&&+2p\int_0^t|Z^\e(s)-Z(s)|^{2p-2}\left<Z^\e(s)-Z(s),\left(\s(\hat{X}^\e_s,\sL_{\hat{X}^\e_s})-\s(X^0_s,\d_{X^0_s})\right)\dif W(s)\right>\\
&&+2p(p-1)\int_0^t|Z^\e(s)-Z(s)|^{2p-4}\Bigl<Z^\e(s)-Z(s),\left(\s(\hat{X}^\e_s,\sL_{\hat{X}^\e_s})-\s(X^0_s,\d_{X^0_s})\right)\\
&&\qquad\times \left(\s(\hat{X}^\e_s,\sL_{\hat{X}^\e_s})-\s(X^0_s,\d_{X^0_s})\right)^{*}(Z^\e(s)-Z(s))\Bigr>\dif s\\
&&+p\int_0^t|Z^\e(s)-Z(s)|^{2p-2}\|\s({\hat{X}^\e_s},\sL_{\hat{X}^\e_s})-\s(X^0_s,\d_{X^0_s})\|^2\dif s.
\de
Then Lemma \ref{L1} implies that
\ce
-2p\int_0^t|Z^\e(s)-Z(s)|^{2p-2}\Bigl<Z^\e(s)-Z(s),\dif \hat{K}^\e(s)-\dif \hat{K}^0(s)\Bigr>\leq 0.
\de 
And by simple computation, we have that 
\ce
&&\frac{b(\hat{X}^\e_s,\sL_{\hat{X}^\e_s})-b(X^0_s,\d_{X^0_s})}{\sqrt{\e}}-{_{\sC^{\ast}}}\<Db(X^0_s,\d_{X^0_s}),Z_s\>_{\sC}-\mE[_{\sC^{\ast}}\<D^Lb(X^0_s,\d_{X^0_s})(X^0_s),Z_s\>_{\sC}]\\
&=&\frac{b(\hat{X}^\e_s,\sL_{\hat{X}^\e_s})-b(X^0_s,\sL_{\hat X^\e_s})}{\sqrt{\e}}-{_{\sC^{\ast}}}\<Db(X^0_s,\sL_{\hat{X}^\e_s}),Z^\e_s\>_{\sC}\\
&&+{_{\sC^{*}}}\<Db(X^0_s,\sL_{\hat{X}^\e_s}),Z^\e_s\>_{\sC}-{_{\sC^{*}}}\<Db(X^0_s,\sL_{\hat{X}^\e_s}),Z_s\>_{\sC}\\
&&+{_{\sC^{*}}}\<Db(X^0_s,\sL_{\hat{X}^\e_s}),Z_s\>_{\sC}-{_{\sC^{*}}}\<Db(X^0_s,\d_{X^0_s}),Z_s\>_{\sC}\\
&&+\frac{b(X^0_s,\sL_{\hat{X}^\e_s})-b(X^0_s,\d_{X^0_s})}{\sqrt{\e}}-\mE[_{\sC^{\ast}}\<D^Lb(X^0_s,\d_{X^0_s})(X^0_s),Z^\e_s\>_{\sC}]\\
&&+\mE[_{\sC^{\ast}}\<D^Lb(X^0_s,\d_{X^0_s})(X^0_s),Z^\e_s\>_{\sC}]-\mE[_{\sC^{\ast}}\<D^Lb(X^0_s,\d_{X^0_s})(X^0_s),Z_s\>_{\sC}].
\de
So, it holds that
\ce
&&\mE\(\sup_{s\in[0,t]}|Z^\e(s)-Z(s)|^{2p}\)\\
&\leq&2p\mE\int_0^t|Z^\e(s)-Z(s)|^{2p-1}\Bigl|\frac{b(\hat{X}^\e_s,\sL_{\hat{X}^\e_s})-b(X^0_s,\sL_{\hat{X}^\e_s})}{\sqrt{\e}}-{_{\sC^{\ast}}}\<Db(X^0_s,\sL_{\hat{X}^\e_s}),Z^\e_s\>_{\sC}\Bigr|\dif s\\ 
&&+2p\mE\int_0^t|Z^\e(s)-Z(s)|^{2p-1}\left|{_{\sC^{\ast}}}\<Db(X^0_s,\sL_{\hat{X}^\e_s}),Z^\e_s\>_{\sC}-{_{\sC^{\ast}}}\<Db(X^0_s,\sL_{\hat{X}^\e_s}),Z_s\>_{\sC}\right|\dif s\\
&&+2p\mE\int_0^t|Z^\e(s)-Z(s)|^{2p-1}\left|{_{\sC^{\ast}}}\<Db(X^0_s,\sL_{\hat{X}^\e_s}),Z_s\>_{\sC}-{_{\sC^{\ast}}}\<Db(X^0_s,\d_{X^0_s}),Z_s\>_{\sC}\right|\dif s\\
&&+2p\mE\int_0^t|Z^\e(s)-Z(s)|^{2p-1}\Bigl|\frac{b(X^0_s,\sL_{\hat{X}^\e_s})-b(X^0_s,\d_{X^0_s})}{\sqrt{\e}}-\mE[_{\sC^{\ast}}\<D^Lb(X^0_s,\d_{X^0_s})(X^0_s),Z^\e_s\>_{\sC}]\Bigr|\dif s\\
&&+2p\mE\int_0^t|Z^\e(s)-Z(s)|^{2p-1}\Bigl|\mE[_{\sC^{\ast}}\<D^Lb(X^0_s,\d_{X^0_s})(X^0_s),Z^\e_s\>_{\sC}]-\mE[_{\sC^{\ast}}\<D^Lb(X^0_s,\d_{X^0_s})(X^0_s),Z_s\>_{\sC}]\Bigr|\dif s\\
&&+2p\mE\sup_{s\in[0,t]}\left|\int_0^s|Z^\e(r)-Z(r)|^{2p-2}\Bigl<Z^\e(r)-Z(r),\left(\s(\hat{X}^\e_r,\sL_{\hat{X}^\e_r})-\s(X^0_r,\d_{X^0_r})\right)\dif W(r)\Bigr>\right|\\
&&+p(2p-1)\mE\int_0^t|Z^\e(s)-Z(s)|^{2p-2}\|\s({\hat{X}^\e_s},\sL_{\hat{X}^\e_s})-\s(X^0_s,\d_{X^0_s})\|^2\dif s\\
&=:&I_1+I_2+I_3+I_4+I_5+I_6+I_7.
\de

For $I_1$, by the Young inequality and $({\bf H}_3)$, it holds that
\ce
I_1&\leq&C\mE\int_0^t\left|\int_0^1{_{\sC^{\ast}}}\<Db(r \sqrt{\e}Z^\e_s+X^0_s,\sL_{\hat{X}^\e_s}),Z^\e_s\>_{\sC}\dif r-{_{\sC^{\ast}}}\<Db(X^0_s,\sL_{\hat{X}^\e_s}),Z^\e_s\>_{\sC}\right|^{2p}\dif s\\
&&+C\mE\int_0^t|Z^\e(s)-Z(s)|^{2p}\dif s\\
&\leq&C\mE\int_0^t\left|\int_0^1\|Db(r\sqrt{\e} Z^\e_s+X^0_s,\sL_{\hat{X}^\e_s})-Db(X^0_s,\sL_{\hat{X}^\e_s})\|_{\sC^{*}}\|Z^\e_s\|_{\infty}\dif r\right|^{2p}\dif s\\
&&+C\mE\int_0^t|Z^\e(s)-Z(s)|^{2p}\dif s\\
&\leq&C\mE\int_0^t\left|\int_0^1 L_3^{1/2}r\sqrt{\e} \|Z^\e_s\|_{\infty}\dif r \cdot \|Z^\e_s\|_{\infty}\right|^{2p}\dif s\\
&&+C\mE\int_0^t|Z^\e(s)-Z(s)|^{2p}\dif s\\
&\leq& C\e^p\int_0^t\mE\|Z^\e_s\|^{4p}_{\infty}\dif s+C\int_0^t\mE|Z^\e(s)-Z(s)|^{2p}\dif s.
\de

As for $I_2$, by the Young inequality, $(\ref{dblinegrow})$, Lemma \ref{lemma:clt-X0bnd} and Lemma \ref{lemma:clt-Zbnd}, it holds that
\ce
I_2&\leq&C\mE\int_0^t|Z^\e(s)-Z(s)|^{2p}\dif s+C\mE\int_0^t\|Db(X^0_s,\sL_{\hat{X}^\e_s})\|^{2p}_{\sC^{\ast}}\|Z^\e_s-Z_s\|^{2p}_{\infty}\dif s\\
&\leq&C\mE\int_0^t|Z^\e(s)-Z(s)|^{2p}\dif s+C\mE\int_0^tL_3^p\left(1+\|X^0_s\|^2_{\infty}+\sL_{\hat{X}^\e_s}(\|\cdot\|_{\infty}^{2})\right)^p\|Z^\e_s-Z_s\|^{2p}_{\infty}\dif s\\
&\leq&C\mE\int_0^t|Z^\e(s)-Z(s)|^{2p}\dif s+C\mE\int_0^t\left(1+\|X^0_s\|^{2p}_{\infty}+\mE\|\hat{X}^\e_s\|^{2p}_{\infty}\right)\|Z^\e_s-Z_s\|^{2p}_{\infty}\dif s\\
&\leq&C\mE\int_0^t|Z^\e(s)-Z(s)|^{2p}\dif s+C\sup_{s\in[0,T]}\left(1+(2^{2p-1}+1)\|X^0_s\|^{2p}_{\infty}+2^{2p-1}\e^p\mE\|Z^\e_s\|^{2p}_{\infty}\right)\\
&&\qquad\times\int_0^t\mE\|Z^\e_s-Z_s\|^{2p}_{\infty}\dif s\\
&\leq&C\int_0^t\mE\sup_{u\in[-r_0,s]}|Z^\e(u)-Z(u)|^{2p}\dif s.
\de

For $I_3$, by the Young inequality and $({\bf H}_3)$, we obtain that
\ce
I_3&\leq&C\mE\int_0^t|Z^\e(s)-Z(s)|^{2p}\dif s+C\mE\int_0^t\|Db(X^0_s,\sL_{\hat{X}^\e_s})-Db(X^0_s,\d_{X^0_s})\|^{2p}_{\sC^{\ast}}\|Z_s\|^{2p}_{\infty}\dif s\\
&\leq&C\mE\int_0^t|Z^\e(s)-Z(s)|^{2p}\dif s+C\mE\int_0^t L^p_3\mW_2(\sL_{\hat{X}^\e_s},\d_{X^0_s})^{2p}\|Z_s\|^{2p}_{\infty}\dif s\\
&\leq&C\mE\int_0^t|Z^\e(s)-Z(s)|^{2p}\dif s+C\mE\int_0^tL_3^p\(\mE\|\hat{X}^\e_s-X^0_s\|^{2p}_{\infty}\)\|Z_s\|_{\infty}^{2p}\dif s\\
&\leq&C\mE\int_0^t|Z^\e(s)-Z(s)|^{2p}\dif s+C\e^p\int_0^t\mE\|Z^\e_s\|^{2p}_{\infty}\mE\|Z_s\|^{2p}_{\infty}\dif s\\
&\leq&C\int_0^t\mE\sup_{u\in[-r_0,s]}|Z^\e(u)-Z(u)|^{2p}\dif s+C\e^p\int_0^t\left(\mE\|Z^\e_s\|^{4p}_{\infty}+\mE\|Z_s\|^{4p}_{\infty}\right)\dif s.
\de

For $I_4$, according to the Young inequality and $({\bf H}_4)$, we have
\ce
I_4&\leq&C\mE\int_0^t\left|\frac{b(X^0_s,\sL_{\hat{X}^\e_s})-b(X^0_s,\d_{X^0_s})}{\sqrt{\e}}-\mE[_{\sC^{\ast}}\<D^Lb(X^0_s,\d_{X^0_s})(X^0_s),Z^\e_s\>_{\sC}]\right|^{2p}\dif s\\
&&+C\mE\int_0^t|Z^\e(s)-Z(s)|^{2p}\dif s\\
&\leq&C\mE\int_0^t\left|\int_0^1\mE[_{\sC^{\ast}}\<D^Lb(X^0_s,\sL_{R_s(r)})(R_s(r)),Z^\e_s\>_{\sC}]\dif r-\mE[_{\sC^{\ast}}\<D^Lb(X^0_s,\d_{X^0_s})(X^0_s),Z^\e_s\>_{\sC}]\right|^{2p}\dif s\\
&&+C\mE\int_0^t|Z^\e(s)-Z(s)|^{2p}\dif s\\
&\leq&C\mE\int_0^t\left|\int_0^1 L_4\left(\mW_2(\sL_{R_s(r)},\d_{X^0_s})+(\mE\|R_s(r)-X^0_s\|^2_{\infty})^{1/2}\right)(\mE\|Z^\e_s\|^2_{\infty})^{1/2} \dif r\right|^{2p}\dif s\\
&&+C\mE\int_0^t|Z^\e(s)-Z(s)|^{2p}\dif s\\
&\leq&C\mE\int_0^t\left|\int_0^1 2L_4(\mE\|R_s(r)-X^0_s\|^2_{\infty})^{1/2}(\mE\|Z^\e_s\|^2_{\infty})^{1/2} \dif r\right|^{2p}\dif s\\
&&+C\mE\int_0^t|Z^\e(s)-Z(s)|^{2p}\dif s\\
&\leq&C\int_0^t\left| 2L_4\sqrt{\e}\mE\|Z^\e_s\|^2_{\infty}\right|^{2p}\dif s+C\mE\int_0^t|Z^\e(s)-Z(s)|^{2p}\dif s\\
&\leq&C\e^p\int_0^t\mE\|Z^\e_s\|^{4p}_{\infty}\dif s+C\int_0^t\mE\sup_{u\in[-r_0,s]}|Z^\e(u)-Z(u)|^{2p}\dif s,
\de
where $R_s(r):=X^0_s+r(\hat{X}^\e_s-X^0_s), r\in[0,1]$.

For $I_5$, $({\bf H}_4)$ and the H\"older inequality imply that
\ce
I_5&\leq&2p\mE\int_0^t|Z^\e(s)-Z(s)|^{2p-1}\left|\mE_{\sC^{\ast}}\<D^Lb(X^0_s,\d_{X^0_s})(X^0_s),Z^\e_s-Z_s\>_{\sC}\right|\dif s\\
&\leq&2p\mE\int_0^t|Z^\e(s)-Z(s)|^{2p-1}\mE\|D^Lb(X^0_s,\d_{X^0_s})\|_{T^*_{\d_{X^0_s},2}}\|Z^\e_s-Z_s\|_{\infty}\dif s\\
&\leq& C\int_0^t\mE\sup_{u\in[-r_0,s]}|Z^\e(u)-Z(u)|^{2p}\dif s.
\de

For $I_6$, applying the Burkholder-Davis-Gundy inequality, the Young inequality and the H$\ddot{\rm o}$lder inequality, we have that
\ce
I_6
&\leq&C\mE\left[\int_0^t|Z^\e(s)-Z(s)|^{4p-2}\left\|\s(\hat{X}^\e_s,\sL_{\hat{X}^\e_s})-\s(X^0_s,\d_{X^0_s})\right\|^2\dif s\right]^{1/2}\\
&\leq&C\mE\left[\sup_{s\in[0,t]}|Z^\e(s)-Z(s)|^{p}\left(\int_0^t|Z^\e(s)-Z(s)|^{2p-2}\left\|\s(\hat{X}^\e_s,\sL_{\hat{X}^\e_s})-\s(X^0_s,\d_{X^0_s})\right\|^2\dif s\right)^{1/2}\right]\\
&\leq&\frac{1}{2}\mE\sup_{s\in[0,t]}|Z^\e(s)-Z(s)|^{2p}+C\mE\int_0^t|Z^\e(s)-Z(s)|^{2p-2}\left\|\s(\hat{X}^\e_s,\sL_{\hat{X}^\e_s})-\s(X^0_s,\d_{X^0_s})\right\|^2\dif s\\
&\leq&\frac{1}{2}\mE\sup_{s\in[0,t]}|Z^\e(s)-Z(s)|^{2p}+C\mE\int_0^t|Z^\e(s)-Z(s)|^{2p}\dif s\\
&&+C\mE\int_0^t\left\|\s(\hat{X}^\e_s,\sL_{\hat{X}^\e_s})-\s(X^0_s,\d_{X^0_s})\right\|^{2p}\dif s.
\de
Besides, $({\bf H}^{\prime}_2)$ yields that
\be
\mE\int_0^t\left\|\s(\hat{X}^\e_s,\sL_{\hat{X}^\e_s})-\s(X^0_s,\d_{X^0_s})\right\|^{2p}\dif s&\leq&{L_2^{\prime}}^p\mE\int_0^t\left(\|\hat{X}^\e_s-X^0_s\|^2_{\infty}+\mE\|\hat{X}^\e_s-X^0_s\|_{\infty}^2\right)^{p}\dif s\no\\
&\leq&C\e^p\int_0^t\mE\|Z^\e_s\|^{2p}_{\infty}\dif s.
\label{eqcltsigmadis}
\ee
By (\ref{eqcltsigmadis}), we derive that
\ce
I_6&\leq&\frac{1}{2}\mE\sup_{s\in[0,t]}|Z^\e(s)-Z(s)|^{2p}+C\mE\int_0^t\sup_{u\in[-r_0,s]}|Z^\e(u)-Z(u)|^{2p}\dif s\\
&&+C\e^p\int_0^t\mE\|Z^\e_s\|^{2p}_{\infty}\dif s.
\de

For $I_7$, the Young inequality and (\ref{eqcltsigmadis}) imply that
\ce
I_7&\leq&C\mE\int_0^t|Z^\e(s)-Z(s)|^{2p}\dif s+C\mE\int_0^t\|\s({\hat{X}^\e_s},\sL_{\hat{X}^\e_s})-\s(X^0_s,\d_{X^0_s})\|^{2p}\dif s\\
&\leq&C\int_0^t\mE\sup_{u\in[-r_0,s]}|Z^\e(u)-Z(u)|^{2p}\dif s+C\e^p\int_0^t\mE\|Z^\e_s\|^{2p}_{\infty}\dif s.
\de 

Collecting the above estimates together, we conclude that
\ce
\mE\(\sup_{s\in[0,t]}|Z^\e(s)-Z(s)|^{2p}\)&\leq&\frac{1}{2}\mE\sup_{s\in[0,t]}|Z^\e(s)-Z(s)|^{2p}+C\int_0^t\mE\sup_{u\in[-r_0,s]}|Z^\e(u)-Z(u)|^{2p}\dif s\\
&&+C\e^p\int_0^t\left(\mE\|Z^\e_s\|^{2p}_{\infty}+\mE\|Z^\e_s\|^{4p}_{\infty}+\mE\|Z_s\|^{4p}_{\infty}\right)\dif s.
\de
Note that $Z^\e(t)=Z^0(t)=0,t\in[-r_0,0]$. Thus, Lemma \ref{lemma:clt-Zbnd} and the Gronwall inequality yield that
\ce
\mE\(\sup_{s\in[-r_0,T]}|Z^\e(s)-Z(s)|^{2p}\)\leq C\e^p.
\de
The proof is complete.
\end{proof}

Finally, we give an example to explain our result.

\bx
Assume that $d=m=1$, and a function $\varphi:\mR\to(-\infty, +\infty]$ is lower semicontinuous convex, and ${\rm Int}(Dom(\varphi))\neq \emptyset$. Then by \cite[Example 2.2]{flqz} we know that $\partial \varphi$ is a maximal monotone operator. 

Next, we take $A=\partial \varphi$, require $0\in\overline{\cD(\partial \varphi)}$ and study the following equation:
\ce\left\{\begin{array}{l}
\dif X^{\e}(t)\in -\partial \varphi(X^{\e}(t))\dif t+f\left(X^{\e}_t+\int_\sC\phi(u)\sL_{X^{\e}_t}(\dif u)\right)\dif t+\sqrt{\e}g(X^{\e}_t)\dif W(t), ~~t\in(0,T], \\
 X^{\e}(t)=\xi(t)\in\overline{\cD(\partial \varphi)},\quad t\in[-r_0,0],
\end{array}
\right.
\de
where $f, g: \sC\to \mR$ and $\phi: \sC\to\sC$ are all twice continuously differentiable mappings with bounded derivatives. So, we know that for any $\zeta\in\sC, \mu\in\sP_2^\sC$
\ce
b(\zeta,\mu)=f\left(\zeta+\int_\sC\phi(u)\mu(\dif u)\right), \quad \sigma(\zeta,\mu)=g(\zeta).
\de
Note that
\ce
&&Db(\zeta,\mu)=Df\left(\zeta+\int_\sC\phi(u)\mu(\dif u)\right), \quad D\sigma(\zeta)=Dg(\zeta),\\
&&D^Lb(\zeta,\mu)(\eta)=Df\left(\zeta+\int_\sC\phi(u)\mu(\dif u)\right)D\phi(\eta).
\de
Thus, it is easy to justify that $({\bf H}_1), ({\bf H}^{\prime}_2), ({\bf H}_3)$ and $({\bf H}_4)$ hold. By Theorem \ref{thm:clt}, we conclude that
\ce
\frac{\hat{X}^{\e}(\cdot)-X^0(\cdot)}{\sqrt{\e}}\xrightarrow{d} Z(\cdot),
\de
where $(Z_{\cdot},\hat{K}^0(\cdot))$ solves the following equation
\ce\left\{\begin{array}{ll}
\dif Z(t)\in -\partial \varphi(Z(t))\dif t+{_{\sC^{\ast}}}\<Df(X^0_t+\phi(X^0_t)),Z_t\>_\sC\dif t+\mE[{_{\sC^{\ast}}}\<Df(X^0_t+\phi(X^0_t))D\phi(X^0_t),Z_t\>_{\sC}]\dif t\\
\quad\quad\quad\quad+g(X^0_t)\dif W(t), \quad t\in(0,T],\\
Z(t)=0,\quad t\in[-r_0,0].
\end{array}
\right.
\de

\ex

\end{document}